\documentclass[reqno]{amsart}

\usepackage{latexsym,amssymb,amsmath, amsfonts}
\usepackage{amsthm}
\usepackage{amscd}
\usepackage[shortlabels]{enumitem}
\usepackage{a4wide}
\usepackage[normalem]{ulem}

\usepackage[backref=page]{hyperref}

\usepackage{graphicx}
\usepackage{import}
\usepackage{adjustbox}

\numberwithin{equation}{section}
\usepackage[dvipsnames,svgnames,table]{xcolor}
\usepackage{hyperref}
\hypersetup{colorlinks=true,
	linkcolor=blue,
	urlcolor=blue,
	citecolor=red} 
\usepackage{import}
\usepackage{transparent}

\newcommand{\R}{\mathbb{R}}

\renewcommand{\le}{\leqslant}
\renewcommand{\ge}{\geqslant}
\renewcommand{\leq}{\leqslant}
\renewcommand{\geq}{\geqslant}

\newcommand{\be}{\begin{equation}}
\newcommand{\en}{\end{equation}}
\newcommand{\ee}{\end{equation}}

\newcommand{\bt}{\begin{theorem}}
\newcommand{\et}{\end{theorem}}
\newcommand{\bp}{\begin{proof}}
\newcommand{\ep}{\end{proof}}
\newcommand{\bc}{\begin{cor}}
\newcommand{\ec}{\end{cor}}
\newcommand{\bl}{\begin{lemma}}
\newcommand{\el}{\end{lemma}}
\newcommand{\bprop}{\begin{prop}}
\newcommand{\eprop}{\end{prop}}

\newtheorem{theorem}{Theorem}[section]
\newtheorem{remark}{Remark}
\newtheorem{lemma}[theorem]{Lemma}
\newtheorem{definition}{Definition}
\newtheorem{proposition}[theorem]{Proposition}
\newtheorem{con}{Conjecture}
\newtheorem{example}[con]{Example}
\newtheorem{corollary}[theorem]{Corollary}
\numberwithin{theorem}{section} \numberwithin{definition}{section}

\newcommand{\RNum}[1]{\uppercase\expandafter{\romannumeral #1\relax}}
\def\R{\mathbb{R}}

\newcommand{\vertiii}[1]{{\left\vert\kern-0.25ex\left\vert\kern-0.25ex\left\vert #1 
		\right\vert\kern-0.25ex\right\vert\kern-0.25ex\right\vert}}

\theoremstyle{definition}

\author[J. Angulo]{Jaime Angulo Pava}

\author[A. Mu\~noz]{Alexander Mu\~noz}

\date{}
\title[NLS on looping-edge graphs] {Nonlinear Schr\"odinger Equations on looping-edge graphs with $\delta'$-type interactions }

\begin{document}

\maketitle

\centerline{ Department of Mathematics,
IME-USP}
 \centerline{Rua do Mat\~ao 1010, Cidade Universit\'aria, CEP 05508-090,
 S\~ao Paulo, SP (Brazil)}
 \centerline{\tt angulo@ime.usp.br, alexd@usp.br}

\section*{Abstract}

In this work, we study the existence and orbital (in)stability of certain standing-wave solutions for the cubic nonlinear Schr\"odinger equation (NLS) posed on a looping-edge graph $\mathcal{G}$, consisting of a circle and a finite number $N$ of infinite half-lines attached to a common vertex. We consider the self-adjoint realization $(\mathcal{H}_Z, D(\mathcal{H}_Z))$ of the Laplacian, where the domain $D(\mathcal{H}_Z)$ encodes on the half-lines a $\delta'$-type vertex conditions  (continuity of derivatives at the vertex, without requiring continuity of the wave function) and $Z \in \mathbb{R}\setminus\{0\}$. On the circle, we propose Jacobian elliptic profiles of dnoidal type combined with either trivial (zero) or soliton tail profiles on the half-lines with full derivative matching at the boundary. For the trivial tail case we establish orbital stability for all $Z \in \mathbb{R}\setminus\{0\}$, while for the non-trivial tail case (which requires $Z < 0$) we establish both existence and orbital (in)stability depending on the relative size of $N$, $Z$, and the phase velocity of the standing wave.

\qquad\\
\textbf{Mathematics Subject Classification (2020)}. Primary
35Q51, 35Q55, 81Q35, 35R02; Secondary 47E05.\\
\textbf{Key words}. Schr\"odinger equation, quantum graphs, boundary systems, extension theory of symmetric operators, Morse and nullity indices for Schr\"odinger operator on non-compact graphs.

\section{Introduction}
The analysis of nonlinear evolution PDE models on metric graphs has potential applications in studying particle and wave dynamics in branched structures and networks. Since such structures appear in various areas of contemporary physics---with applications in electronics, biology, materials science, and nanotechnology---the development of effective modeling tools is essential for addressing the many practical problems that arise in these fields (see \cite{BK} and references therein).

However, real systems can exhibit strong inhomogeneities due to varying nonlinear coefficients across different regions of the spatial domain or due to the specific geometry of the domain itself. To explore key features of the nonlinear Schr\"odinger equation (NLS model), we will therefore choose a `simple" metric graph tool, such as the looping-edge graph (see Figure 1). The flexibility in designing the graph's geometry enables the investigation of a wide range of dynamic behaviors.

Recently, nonlinear models on graphs, such as the nonlinear Schr\"odinger equation, the sine-Gordon model, and the Korteweg-de Vries model, have been studied extensively (see \cite{AdaNoj14, AngGol17a, AngGol17b, AC, AC1, AP1, AP2, AP3, AST, Fid15, Noj14} and references therein). From a mathematical viewpoint, an evolution model on a graph is equivalent to a system of PDEs defined on the edges (intervals), where the coupling is determined exclusively by the boundary conditions at the vertices (known as the "topology of the graph"), which governs the dynamics on the network via groups. This area of study has attracted significant attention, particularly in the context of soliton transport. Solitons and other nonlinear waves in branched systems provide valuable insights into the dynamics of these models. To mention a few examples of bifurcated systems: Josephson junction structures, networks of planar optical waveguides and fibers, branched structures associated with DNA, blood pressure waves in large arteries, nerve impulses in complex arrays of neurons, conducting polymers, and Bose-Einstein condensation (see  \cite{BK, BlaExn08, BurCas01, Chuiko, Crepeau, Ex, Fid15, K, Mug15, Noj14, SBM} and reference therein).

Addressing these issues is challenging because both the equations of motion and the graph topology can be complex. Additionally, a critical aspect that complicates the analysis is the presence of one or more vertices where a soliton profile, arriving at the vertex along one of the edges, undergoes complicated behaviors such as reflection and the emergence of radiation. This makes it difficult to observe how energy travels across the network. As a result, the study of soliton propagation through networks presents significant challenges. The mechanisms for the existence and stability (or instability) of soliton profiles remain unclear for many types of graphs and models.

We recall that a metric graph $\mathbb{G}$ is a structure represented by a finite number of vertices $V=\{\nu_i\}$ and a set of adjacent edges at the vertices $E=\{e_j\}$ (for further details, see \cite{BK}). Each edge $e_j$ can be identified with a finite or infinite interval of the real line, $I_e$. Thus, the edges of $\mathbb{G}$ are not merely abstract relations between vertices but can be viewed as physical "wires" or "networks" connecting them. The notation $e \in E$ will be used to indicate that $e$ is an edge of $\mathbb{G}$. This identification introduces the coordinate $x_e$ along the edge $e$.

\begin{figure}
    \centering
    \includegraphics[width=0.75\linewidth]{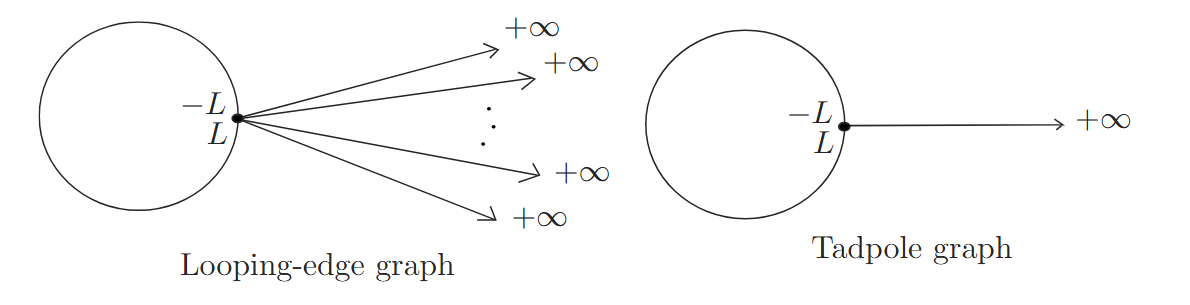}
    \caption{}
    \label{Fig1}
\end{figure}

On this work we will focus mainly on  looping-edge graphs composed of a circle and a finite number $N$ of infinite half-lines attached to a common vertex. If the circle edge is identified with the interval $[-L, L]$ and the half-lines with $[L, \infty)$, we obtain a particular metric graph structure, that we will denote with $\mathcal{G}$, represented by $\mathcal{G}=V\cup E$ where $E = \{[-L, L], [L, \infty), \dots, [L, \infty)\}$, $V = \{\nu\}$ and along the coordinate systems of each  edge $e$, the single vertex $\nu$ is located exactly at $x_e=L$. (see Figure \ref{Fig1}). We use the subscripts $0$ and $j$ to refer to the edges $e_0 = [-L, L]$ and $e_j = [L, \infty)$, $j = 1, \dots, N$, respectively. The special case $N = 1$ is known as a tadpole graph (see also Figure \ref{Fig1}).

A wave function $\textbf{U}$ defined on $\mathcal{G}$ will be understood as an $(N+1)$-tuple of functions $\textbf{U} = (\phi, (\psi_j)_{j=1}^N) = (\phi, \psi_1, \dots, \psi_N)$, where $\phi$ is defined on $e_0 = [-L, L]$ and $\psi_j$ is defined on $e_j = [L, \infty)$ for $j = 1, \dots, N$. 

We discuss in this work the case of nonlinear Schr\"odinger-type equations (NLS)
\begin{equation}\label{NLS}
i\textbf{U}_t\ + \Delta\textbf{U}+ |\textbf{U}|^{p-1}\textbf{U} =\textbf{0},\quad p>1,
\end{equation}
where the action of the Laplacian operator $\Delta$ on a general graph $\mathbb{G}$ is given by
\begin{equation}
-\Delta: (u_e)_{e\in E} \to (-u''_e)_{e\in E}.
\end{equation}
The NLS in \eqref{NLS} has been extensively studied in recent literature for various types of metric graphs $\mathbb{G}$ and with specific domains for the Laplacian, which make it a self-adjoint operator on $L^2(\mathbb{G})$ (see the review manuscript \cite{KNP}). For instance, research has focused on star graphs (\cite{AdaNoj15, AdaNoj14, AngGol17a, AngGol17b, KPG, Noj14} and references therein), flower graphs, looping-edge graphs, dumbbell graphs, double-bridge graphs, and periodic ring graphs (\cite{An1, An2, BMP, CFN, KP, KMPX, NPS, NP, Pan} and references therein).

This work contains a characterization of self-adjoint extensions of the symmetric operator $\mathcal H_0\equiv -\Delta$ on the metric graph $\mathcal{G}$ with domain $D(\mathcal H_0)$ defined by
\begin{equation}
    \label{cinf} D(\mathcal H_0)=C_0^\infty(-L,L)\oplus \bigoplus_{j=1}^N C_0^\infty(L,+\infty).
\end{equation}
Since the deficiency indices of $(\mathcal H_0, D(\mathcal H_0))$ are equal to $2+N$, the extension theory of Krein and von Neumann provides a $(2+N)^2$-parameter family of self-adjoint extensions $(-\Delta_{ext}, D(-\Delta_{ext}))$ of $\mathcal{H}_0$ on $L^2(\mathcal{G})$, where the action of $-\Delta_{ext}$ is given by $-\Delta_{ext} = -\Delta$. Each member of this family defines a unitary dynamic for the linear evolution problem
\begin{equation}\label{evol}
\left\{ \begin{array}{ll}
i\textbf U_t = -\Delta_{ext} \textbf U\\
\textbf U(0)= \textbf u_0\in D(-\Delta_{ext}).
\end{array} \right.
\end{equation}

The nature of the self-adjoint extensions governing the dynamics of \eqref{evol} may vary significantly. If both the looping-edge geometry $\phi(L)=\phi(-L)$ and the continuity conditions at the vertex are preserved, we obtain a one-parameter family of self-adjoint extensions $(-\Delta_{\mathrm{ext}}, D_{Z,N})_{Z \in \mathbb{R}}$ given by
\begin{equation}\label{DomainN0}
\begin{aligned}
D_{Z, N} = \bigg\{\mathbf{U} \in H^2(\mathcal{G}) \ \bigg|\ &
\phi(L)=\phi(-L)=\psi_1(L)=\cdots=\psi_N(L), \\
& \phi'(L)-\phi'(-L)=\sum_{j=1}^N \psi_j'(L)+Z\,\psi_1(L) \bigg\}.
\end{aligned}
\end{equation}

Several studies related to the NLS with the domain $D_{Z, N}$ have been carried out, particularly concerning the existence and stability of standing-wave solutions (see \cite{AST, ASTcri, ASTmul, An1, An2, Ardila, CFN, CFN2, KNP, NP, NPS}, and references therein). The boundary conditions in \eqref{DomainN0} are referred to as \emph{$\delta$-type} if $Z \neq 0$, and as \emph{Neumann--Kirchhoff type} if $Z = 0$.

Domains as in \eqref{DomainN0} may appear in conceivable experiments of electron motion in thin metallic structures, so-called ``network of quantum wires'' (see \cite{Ex}, \cite{ExSere}, and references therein).

Our focus in this manuscript is on the existence and orbital (in)stability of standing-wave solutions to the cubic NLS model on looping-edge graphs with derivative coupling interactions at the vertex. As a first step in the analysis of such interactions on looping-edge graphs, we consider the Laplacian as a self-adjoint extension with domain \begin{equation}\label{eq:DZ-final_intro}\begin{split}
    D(\mathcal{H}_Z)
  = \Bigl\{\mathbf{U}\in H^2(\mathcal{G})\;\Big|\;
        \phi(-L)&=\phi(L),\quad
        \phi'(-L)=\phi'(L),\quad\\&
        \psi_1'(L)=\cdots=\psi_N'(L),\quad
        \sum_{j=1}^N\psi_j(L)=Z\,\psi_1'(L)
    \Bigr\}.
\end{split}
\end{equation}
obtained in \cite{AM1} via the theory of boundary systems \cite{Schu2015}.

We study standing-wave solutions are of the form
\[
\mathbf{U}(x,t)=e^{i\omega t}\Theta(x),
\]
with $\omega>0$, where $\Theta=(\Phi,\Psi_1,\Psi_2,\ldots,\Psi_N)$ has real-valued components, $\Phi:[-L,L]\to\mathbb{R}$ and $\Psi_i:[L,+\infty)\to\mathbb{R}$, and satisfies the stationary NLS vectorial equation
\begin{equation}\label{standingN0}
 -\Delta \Theta + \omega \Theta - |\Theta|^{2}\Theta = \mathbf{0},
\end{equation}
with $\Theta \in D(\mathcal{H}_Z)$ (see \eqref{eq:DZ-final_intro}).

\begin{figure} 
\centering
    \includegraphics[width=0.5\linewidth]{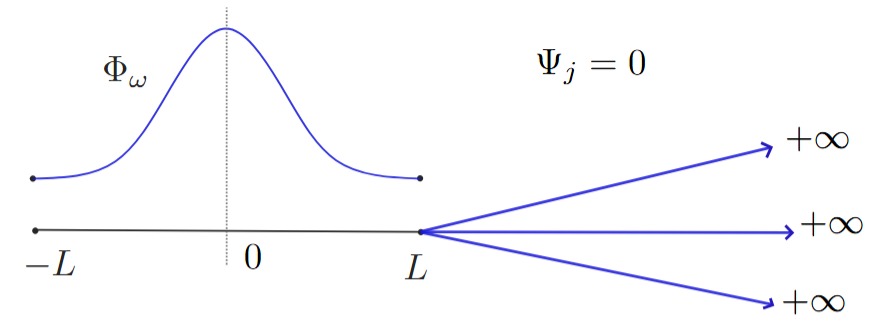}
    \caption{$\Phi_\omega$ dnoidal $+$ trivial soliton profiles.}
    \label{Fig2i}
\end{figure}

\begin{figure}
    \centering
    \includegraphics[width=0.5\linewidth]{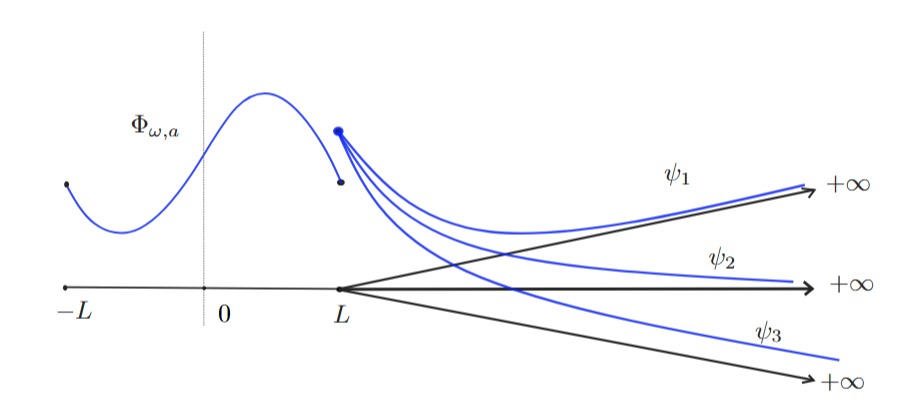}
    \caption{$\Phi_\omega$-dnoidal translated $+$ tail profiles.}
    \label{Fig3i}
\end{figure}

In this work, we consider particular profiles of the form $(\Phi_1,\mathbf{0})$ and $(\Phi_2, (\Psi)_{i=1}^N)$, where $\Phi_i$ are profiles given by Jacobi elliptic functions of \emph{dnoidal} type and $\Psi$ is a soliton tail profile (see Figures~\ref{Fig2i} and~\ref{Fig3i}, respectively). In the case of trivial (zero) solitons on the half-lines, our main existence and stability result can be stated as follows:

\begin{theorem}\label{1stability1}
Let $L>0$ be arbitrary but fixed and $Z\in \mathbb{R}\setminus\{0\}$. For every $\omega>\frac{\pi^2}{2L^2}$ (with $\omega> \frac{N^2}{Z^2}$ if $Z<0$), the standing-wave solution $e^{i\omega t}(\Phi_\omega, \mathbf{0})$ is orbitally stable in $H^1(\mathcal{G})$ under the flow of the cubic NLS equation on a looping-edge graph. Here, $\Phi_\omega$ represents the dnoidal profile
\begin{equation}\label{dnoidal0}
\Phi_\omega(x)=\eta_1\,\mathrm{dn}\Big(\frac{\eta_1}{\sqrt{2}}\, x;k\Big),
\end{equation}
with elliptic modulus $k\in (0,1)$ determined by
\begin{equation}
k^2=\frac{\eta_1^2-\eta_2^2}{\eta_1^2},\qquad \eta_1^2+\eta_2^2=2\omega,\qquad 0<\eta_2<\eta_1.
\end{equation}
\end{theorem}

The proof of Theorem~\ref{1stability1} is divided into Corollary~\ref{d3} and Theorem~\ref{1stability1-ext} along the manuscript.

\medskip
When considering tail profiles on the half-lines, we establish the following result.

\begin{theorem}\label{2stability2}
Let $L>0$, $Z<0$, and $N\geq 1$. For $\omega$ in the open set $I(Z,N,L)$ defined in \eqref{Interval} and $a=a(\omega)\in (0,L)$, there exists a smooth curve of standing-wave solutions for the cubic NLS,
\[
\omega \longmapsto \big(\Phi_{\omega, a}, (\Psi_{\omega, Z, N})_{i=1}^N\big) \in D(\mathcal{H}_Z).
\]
Here, $\Phi_{\omega, a}$ is a dnoidal profile given by $\Phi_{\omega, a}(x)=\Phi_{\omega}(x-a)$ for $x\in [-L,L]$, and for $\omega> \frac{N^2}{Z^2}$
\begin{equation}\label{Ntails0}
\Psi_{\omega, Z, N}(x)=\sqrt{2\omega}\,\mathrm{sech}\Big( \sqrt{\omega}(x-L)+\tanh^{-1}\Big(\frac{-N}{Z\sqrt{\omega}}\Big)\Big),\qquad x\geq L.
\end{equation}

Moreover the standing wave $e^{i\omega t}\big(\Phi_{\omega, a},  (\Psi_{\omega, Z, N})_{i=1}^N\big)$ is orbitally unstable in $H^1(\mathcal{G})$.
\end{theorem}
The proof of Theorem~\ref{2stability2} in presented in the manuscript in two steps related to  Theorems~\ref{existence1} and
\ref{thm:instability-HZ}.

\medskip
The statements of the orbital (in)stability results in Theorems~\ref{1stability1} and~\ref{2stability2} are based on the abstract framework established by Grillakis, Shatah, and Strauss in \cite{GSS2, GSS1} adapted to the case of metric graphs. For the convenience of the reader, we present in Theorem~\ref{2main} an adaptation of these abstract results to the particular case of looping-edge graphs. 

In order to apply the theory in \cite{GSS2, GSS1}, it is necessary to establish local and global well-posedness, defined as the existence, uniqueness, and continuous dependence of solutions within the energy space (see \cite{Caz}). In our setting, global well-posedness can be obtained by employing the standard strategy of using the conserved quantities
\begin{equation}\label{Ener}
E_Z(\mathbf{U}) = \frac{1}{2} \|\nabla \mathbf{U}\|^2_{L^2(\mathcal{G})} - \frac{1}{p+1} \|\mathbf{U}\|^{p+1}_{L^{p+1}(\mathcal{G})} + \frac{1}{2Z} \left| \sum_{i=1}^N \psi_i(L) \right|^2,
\qquad Z \neq 0 \quad (\text{energy}),
\end{equation}
and
\begin{equation}\label{mass}
Q(\mathbf{U}) = \|\mathbf{U}\|^2_{L^2(\mathcal{G})}, \hskip2.7in (\text{mass}),
\end{equation}
both defined in the energy space $H^1(\mathcal{G})$, to extend the solution globally in time. We observe that $E_Z$ and $Q$ belong to $C^2(H^1(\mathcal{G}))$ since $p>1$.

We now state the well-posedness result for \eqref{NLS}.

\begin{theorem}[Local and global well-posedness]\label{global}
Let $\mathbf{U}_0 \in H^1(\mathcal{G})$ and $p>1$. Then the integral equation associated with \eqref{NLS},
\begin{equation}\label{inteq}
\mathbf{U}(t) = e^{it\mathcal{H}_Z}\mathbf{U}_0
- i \int_0^t e^{i(t-s)\mathcal{H}_Z } |\mathbf{U}(s)|^{p-1}\mathbf{U}(s)\, ds,
\end{equation}
has a unique solution $\mathbf{U} \in C((-T,T); H^1(\mathcal{G}))$ for some $T = T(\|\mathbf{U}_0\|_{H^1}) > 0$ with initial condition $\mathbf{U}(0)=\mathbf{U}_0$. Moreover, for any $T_0 < T$, the mapping
\[
\mathbf{U}_0 \in H^1(\mathcal{G}) \ \longmapsto\ \mathbf{U} \in C([-T_0,T_0]; H^1(\mathcal{G}))
\]
is continuous. In particular, for $p>2$ this mapping is of class $C^2$. For $1<p<5$, the solution can be extended globally in time, i.e., $T$ can be taken arbitrarily large.
\end{theorem}

For $\mathcal{H}_Z$ denoting the operator $-\Delta$ acting on $L^2(\mathcal{G})$ with domain $D(\mathcal{H}_Z)$, Theorem~\ref{group-ext} in Section~4 provides a representation of the unitary group $e^{it\mathcal{H}_Z}$ to be used in our stability analysis. Consequently, a spectral description of the one-parameter family of self-adjoint operators $(\mathcal{H}_Z, D(\mathcal{H}_Z))_{Z\in \mathbb{R}\setminus\{0\}}$ is required (see Proposition~\ref{pointesp-ext}).

\medskip
\textbf{Organization of the article.} 
Section~2 fixes notation and recalls standard function spaces on metric graphs.
Section~3 develops the boundary-systems framework for the Schr\"odinger operator on $\mathcal{G}$ and properly defines $(\mathcal{H}_Z,D(\mathcal{H}_Z))$. Section~4 is devoted to the well-posedness theory.
Section~5 establishes the existence of the standing-wave profiles on $\mathcal{G}$.
Section~6 contains the orbital (in)stability analysis: the trivial-tail case is treated in Section~6.1 while the non-trivial tail case is treated in Section~6.2.

\section{Notation}
Let $-\infty\leq a<b\leq\infty$. We denote by $L^2(a,b)$  the  Hilbert space equipped with the inner product $(u,v)=\int\limits_a^b u(x)\overline{v(x)}dx$.  By $H^n(\Omega)$  we denote the classical  Sobolev spaces on $\Omega\subset \mathbb R$ with the usual norm. For a  metric looping-edge graph $\mathcal{G}$, the graph $\dot{\mathcal{G}}$ denotes the star graph obtained by removing the circular edge $e_0$ from the looping-edge graph.
On $\mathcal{G}$ we define the  $ L^ p(\mathcal G)$-spaces by $
 L^ p(\mathcal G)=L^ p(-L, L)\oplus \bigoplus_j L^p( L, +\infty)$, $p>1$, with the natural norms. For $\textbf{U}=(f, \{g_j\}), \textbf{V}=(\tilde{f}, \{\tilde{g}_j\}) \in L^2(\mathcal{G})$,  the natural inner product in $L^2(\mathcal{G})$  is defined by
 $[\textbf{U}, \textbf{V}]= \int_{-L}^L f(x)\overline{\tilde{f}(x)}dx + \sum_{j=1}^N \int_L^{\infty} g_j(x)\overline{\tilde{g}_j(x)}dx$.  For any $n\geqq 0$, we have the  $H^n(\mathcal{G})$-Sobolev spaces,
 $$
 H^n(\mathcal{G})=H^ n(-L, L)\oplus \bigoplus_{j=1}^ N H^ n( L, +\infty).
 $$
For $\mathbf{U}=(\phi,\psi_1,\dots,\psi_N)$ posed on $\mathcal{G}$ we understand any derivative $\phi^{(k)}(L)$ as a inward derivative, $\phi^{(k)}(-L)$ as an outward derivative and $\psi_j^{(k)}(L)$ as an outward derivative.
 For $\mathbf{U} = (\phi,\psi_1,\dots,\psi_N)$, we denote by $\vec{\mathbf{U}}$ and $\vec{\mathbf{U}}'$ the row vectors collecting the values of $\mathbf{U}$ and its derivatives at the vertex, respectively:
\begin{equation*}
\begin{split}
    \vec{\mathbf{U}} := \vec{\mathbf{U}}(L) 
= (\phi(-L),\ \phi(L),\ \psi_1(L),\ \dots, \psi_N(L))^T,\\ 
\vec{\mathbf{U}}' := \vec{\mathbf{U}}'(L):=(\phi'(-L),\ \phi'(L),\ \psi_1'(L),\ \dots, \psi_N'(L) )^T.
\end{split}
\end{equation*}

 Let $A$ be a  closed densely defined symmetric operator in the Hilbert space $H$. The domain of $A$ is denoted by $D(A)$. The deficiency indices of $A$ are denoted by  $n_\pm(A):=\dim ker(A^*\mp iI)$, with $A^*$ denoting the adjoint operator of $A$.

 \section{Boundary systems and unitary dynamics for the linear Schr\"odinger model on  looping-edge graphs}\label{S3}

In recent years, the study of the NLS on graphs has made considerable progress, and the development of new insights to understand the dynamics of this model has been very fruitful. It is well known that the boundary conditions at the vertices govern the dynamics on the network. A classical example is the so-called Neumann–Kirchhoff condition (widely used in real-world applications, see for instance \cite{Ex, ExS, ExSere}). In the case of a looping-edge graph, the $\delta$-type boundary conditions (continuity at the vertex) in \eqref{DomainN0} have been studied extensively. 

In this manuscript, our focus is to study certain profiles with derivative matching at the vertex, inspired by the  $\delta'$-type interactions (lack of continuity at the vertex but continuity of derivatives) as in the case of metric star graphs \cite{AngGol17b}. Based on boundary systems, we provide an answer to the question of which are all the self-adjoint extensions of the Laplace operator that include continuity of derivatives at the vertex and are posed on a looping-edge graph. Our approach is based on the abstract extension theory of symmetric operators developed in Section~2 of \cite{Schu2015}. Indeed in \cite{AM1} this abstract extension theory was applied in the study of the Airy and Schr\"odinger operators posed on looping-edge graphs, for sake of completeness we include some details here.

\subsection{Deficiency indices of the free Schr\"odinger operator}
For the symmetric operator $\mathcal{H}_0$ we have that free solutions of $-\partial_x^2 u_e \pm iu_e=0$ are given as linear combinations of complex exponentials of the form $$u_e=c_1 e^{\frac{(1\pm i)}{\sqrt{2}}x_e}+c_2 e^{-\frac{(1\pm i)}{\sqrt{2}}x_e}.$$ 

We note that $e^{\pm\frac{(1\pm i)}{\sqrt{2}}x}$ belong to $L^2(-L,L)$. On the other hand, on a half-line $[L,\infty)$, it would be square integrable only if the real part of $\pm\frac{(1\pm i)}{\sqrt{2}}$ is negative. Along each edge $e=[L,\infty)$ we would therefore always have that the deficiency indices are ${d_+}_e={d_-}_e=1$.

In a looping-edge graph with $N$ half-lines attached at the vertex we would have for $\mathcal{H}_0$ that the deficiency indices  $d_{\mp}=dim (Ker(\mathcal{H}_0\pm iI))$ satisfy $d_+=d_-=2+N$. 
From the von Neumann and Krein theory, one can always expect self-adjoint extensions of $\mathcal{H}_0$ and therefore unitary dynamics.

\subsection{Boundary systems}\label{sec:bdry_schr}
In this subsection we adapt the ideas developed in the case of the Airy operator, for the Schr\"odinger operator posed on $\mathcal{G}$.  

 \begin{definition}
  We define on the graph of $\mathcal{H}_0^*$, $G(\mathcal{H}_0^*)=\{(\textbf{U},\mathcal{H}_0^* \textbf{U})\mid \textbf{U} \in D(\mathcal{H}_0^*)\}$, the standard skew-symmetric form $\Lambda: G(\mathcal{H}_0^*)\times G(\mathcal{H}_0^*) \to \mathbb{C}, $ defined by \begin{equation}\label{form-}
\Lambda((\textbf{U},\mathcal{H}_0^*\textbf{U}),(\textbf{V},\mathcal{H}_0^*\textbf{V})):=[\mathcal{H}_0^*\textbf{U}, \textbf{V}]-[\textbf{U}, \mathcal{H}_0^*\textbf{V}]. 
\end{equation} 
\end{definition}

\begin{definition}\label{def:omega}
   Let $\mathcal{J}=\mathbb{C}^{N+2}$. Define the standard skew-symmetric sesquilinear form on
    $\mathcal{J}\oplus\mathcal{J}$ by
    \begin{equation}\label{eq:omega}
        \rho\bigl((x_1,y_1),(x_2,y_2)\bigr)
        :=\langle y_1,x_2\rangle_{\mathcal{J}}
         -\langle x_1,y_2\rangle_{\mathcal{J}},
        \qquad (x_i,y_i)\in\mathcal{J}\oplus\mathcal{J},
    \end{equation}
    where $\langle\,\cdot\,,\,\cdot\,\rangle_{\mathcal{J}}$ is the
    standard inner product on $\mathbb{C}^{N+2}$.
\end{definition}

\begin{definition}
    We define $\mathbb{F}$ to be the surjective boundary map
    \begin{equation}
        \begin{split}
            \mathbb{F}:&\ G(\mathcal{H}_0^*)\to \mathcal{J}\oplus \mathcal{J}\\
            &(\mathbf{U},\mathcal{H}_0^*\mathbf{U})\mapsto \bigl(\Gamma_0\mathbf{U},\,\Gamma_1\mathbf{U}\bigr),
        \end{split}
    \end{equation}
    where 
    \begin{equation*}
        \Gamma_0\mathbf{U}=\begin{pmatrix}\phi(-L)\\ \phi(L)\\ \psi_1(L)\\ \vdots\\ \psi_N(L)\end{pmatrix},
        \qquad
        \Gamma_1\mathbf{U}=\begin{pmatrix}\phi'(-L)\\ -\phi'(L)\\ \psi_1'(L)\\ \vdots\\ \psi_N'(L)\end{pmatrix}.
    \end{equation*}
\end{definition}

\begin{proposition}\label{prop:boundary_system}
    The quintuple $(\Lambda,\mathcal{J},\mathcal{J},\mathbb{F},\rho)$ is a \emph{boundary system} for $\mathcal{H}_0$, that is,
    $$\Lambda((\textbf{U},\mathcal{H}_0^*\textbf{U}),(\textbf{V},\mathcal{H}_0^*\textbf{V
}))
                =\rho\bigl(\mathbb{F}(\mathbf{U},\mathcal{H}_0^*\mathbf{U}),\,
                         \mathbb{F}(\mathbf{V},\mathcal{H}_0^*\mathbf{V})\bigr).$$
\end{proposition}

\begin{proof}
    The proof follows by integration by parts.
\end{proof}

\subsection{Self-adjoint extensions via boundary systems}

The following theorem, taken from \cite[Theorem~3.3]{Schu2015}, gives a complete
parametrization of all self-adjoint extensions of $\mathcal{H}_0$ in terms of
the boundary system constructed in the previous subsection.

\begin{theorem}\label{thm:SAext}
    An operator $\mathcal{H}$ is a self-adjoint extension of $\mathcal{H}_0$ if
    and only if there exist a subspace $X\subset\mathcal{J}=\mathbb{C}^{N+2}$
    and a self-adjoint operator $\mathcal{L}:X\to X$ such that
    $\mathcal{H}=\mathcal{H}_{X,\mathcal{L}}$, where
    \begin{equation}\label{eq:SAext}
        D(\mathcal{H}_{X,\mathcal{L}})=\Bigl\{\mathbf{U}\in H^2(\mathcal{G})\;\Big|\;
            \Gamma_0\mathbf{U}\in X, \quad \mathcal{L}\Gamma_0\mathbf{U}=P_X\Gamma_1\mathbf{U}
            \Bigr\}.
    \end{equation} 
    Here $P_X:\mathcal{J}\to X$ is the orthogonal projection onto $X$.
\end{theorem}

%

\begin{example}[Self-adjoint extension $\mathcal{H}_Z$ on a
looping-edge graph]
\label{ex:HZ-derivation}

Consider the boundary parameter matrix
\begin{equation}\label{eq:M-choice}
  M = \begin{pmatrix} 0 & \mathbf{0}^{*} \\ \mathbf{0} &
  \dfrac{1}{Z}\,\mathbf{1}\mathbf{1}^{\top} \end{pmatrix},
  \qquad Z \in \mathbb{R}\setminus\{0\},
\end{equation}
where $\mathbf{1}=(1,\ldots,1)^{\top}\in\mathbb{R}^N$.

The looping condition $\phi(-L)=\phi(L)$ constrains the boundary
value $\Gamma_0\mathbf{U}\in\mathcal{J}=\mathbb{C}^{N+2}$ to lie in
the subspace
\[
  X := \bigl\{x\in\mathcal{J}\;\big|\;x_1=x_2\bigr\},
  \qquad
  X^{\perp} = \operatorname{span}\{e_1-e_2\},
  \qquad
  \dim X = N+1.
\]
An orthonormal basis for $X$ is
\begin{equation}\label{eq:basis}
  g_1 := \frac{e_1+e_2}{\sqrt{2}},
  \qquad
  g_k := e_{k+1}, \quad k = 2,\ldots,N+1,
\end{equation}
where $\{e_j\}_{j=1}^{N+2}$ denotes the standard basis of
$\mathbb{C}^{N+2}$.

\medskip
A direct computation shows that the coordinate representations of
$\Gamma_0\mathbf{U}$ and of the projected boundary map
$P_X\Gamma_1\mathbf{U}$ in the basis $\{g_k\}$ are
\begin{equation}\label{eq:tilde-maps}
  \widetilde{\Gamma}_0\mathbf{U}
  = \begin{pmatrix}
      \sqrt{2}\,\phi(-L) \\ \psi_1(L) \\ \vdots \\ \psi_N(L)
    \end{pmatrix},
  \qquad
  \widetilde{\Gamma}_1\mathbf{U}
  = \begin{pmatrix}
      \dfrac{\phi'(-L)-\phi'(L)}{\sqrt{2}} \\[6pt]
      \psi_1'(L) \\ \vdots \\ \psi_N'(L)
    \end{pmatrix}.
\end{equation}

\medskip
The self-adjoint extension corresponding to the parameter matrix $M$
is characterized by the equation
\begin{equation}\label{eq:sa-eq}
  M\,\widetilde{\Gamma}_0\mathbf{U}
  = \widetilde{\Gamma}_1\mathbf{U}.
\end{equation}
For the matrix \eqref{eq:M-choice} this reads
\begin{equation}\label{eq:sa-expanded}
  \begin{pmatrix}
    0 & \mathbf{0}^{*} \\
    \mathbf{0} & \dfrac{1}{Z}\mathbf{1}\mathbf{1}^{\top}
  \end{pmatrix}
  \begin{pmatrix}
    \sqrt{2}\,\phi(-L) \\ \psi_1(L) \\ \vdots \\ \psi_N(L)
  \end{pmatrix}
  =
  \begin{pmatrix}
    \dfrac{\phi'(-L)-\phi'(L)}{\sqrt{2}} \\[6pt]
    \psi_1'(L) \\ \vdots \\ \psi_N'(L)
  \end{pmatrix}.
\end{equation}

\medskip
The $(1,1)$-entry of $M$ is $0$ and the first row of the off-diagonal
block is $\mathbf{0}^*=(0,\ldots,0)$, so the left-hand side of the
first equation is
\[
  0\cdot\sqrt{2}\,\phi(-L) + \mathbf{0}^*\cdot\vec\psi(L) = 0.
\]
Setting this equal to $\tfrac{\phi'(-L)-\phi'(L)}{\sqrt{2}}$ gives
\begin{equation}\label{eq:bc-loop-Z}
  \frac{\phi'(-L)-\phi'(L)}{\sqrt{2}} = 0
  \qquad\Longleftrightarrow\qquad
  \phi'(-L) = \phi'(L).
\end{equation}
Together with the looping condition already imposed, this is the
\emph{full periodic boundary condition} on the loop.

\medskip
For $j=1,\ldots,N$, the $(j+1)$-th row of $M$ is the
$(j+1)$-th row of the block $\tfrac{1}{Z}\mathbf{1}\mathbf{1}^{\top}$,
which equals $\tfrac{1}{Z}(1,1,\ldots,1)\in\mathbb{R}^N$,
and the off-diagonal entry is $(\mathbf{0})_j=0$.
Hence the left-hand side of the $(j+1)$-th equation is
\[
  0\cdot\sqrt{2}\,\phi(-L)
  + \frac{1}{Z}\bigl(1,1,\ldots,1\bigr)\vec\psi(L)
  = \frac{1}{Z}\sum_{i=1}^N \psi_i(L).
\]
Setting this equal to $\psi_j'(L)$ yields, for every
$j=1,\ldots,N$,
\begin{equation}\label{eq:bc-halfline-j}
  \psi_j'(L) = \frac{1}{Z}\sum_{i=1}^N \psi_i(L).
\end{equation}

Since the right-hand side of \eqref{eq:bc-halfline-j} is the
\emph{same} for every $j$, equations \eqref{eq:bc-halfline-j}
are equivalent to the two conditions
\begin{equation}\label{eq:sum_cond-der}
  \psi_1'(L) = \psi_2'(L) = \cdots = \psi_N'(L),
  \qquad
  \sum_{j=1}^N \psi_j(L) = Z\,\psi_1'(L).
\end{equation}

\medskip
Collecting \eqref{eq:bc-loop-Z} and \eqref{eq:sum_cond-der},
and recalling that the looping condition $\phi(-L)=\phi(L)$ was
imposed from the outset, we obtain the self-adjoint extension  $(\mathcal{H}_Z,D(\mathcal{H}_Z))$ with domain
\begin{equation}\label{eq:DZ-final}
\begin{split}
    D(\mathcal{H}_Z)
  = \Bigl\{\mathbf{U}\in H^2(\mathcal{G})\;\Big|\;
        \phi(-L)&=\phi(L),\quad
        \phi'(-L)=\phi'(L),\quad\\&
        \psi_1'(L)=\cdots=\psi_N'(L),\quad
        \sum_{j=1}^N\psi_j(L)=Z\,\psi_1'(L)
    \Bigr\}.
\end{split}
\end{equation}
\end{example}

\section{Global well-posedness of Cauchy Problem}

In this section we study the well-posedness of equation \eqref{NLS} for $p>1$. Before doing so, we first describe the spectrum of the one parameter family of self-adjoint operators $$(-\Delta, D(\mathcal{H}_Z))\qquad \mbox{for }\ Z \in \R\setminus\{0\},\ N\geqq 1,$$ and $D(\mathcal{H}_Z)$ defined in \eqref{eq:DZ-final}.

\begin{proposition}\label{pointesp-ext}
Let $Z \in \mathbb{R}\setminus\{0\}$ and let $\mathcal{H}_Z$ be the
self-adjoint operator $-d^2/dx^2$ with domain $D(\mathcal{H}_Z)$
given by \eqref{eq:DZ-final}.
Then:

\begin{enumerate}
\item[1)] For $Z >0$, the spectrum of $\mathcal{H}_Z$ in
  $L^2(\mathcal{G})$ is
  \[
    \sigma(\mathcal{H}_Z)
    = [0,+\infty)
    = \sigma_c(\mathcal{H}_Z)\cup\sigma_p(\mathcal{H}_Z),
  \]
  where
  \[
    \sigma_p(\mathcal{H}_Z)
    = \Bigl\{\tfrac{n^2\pi^2}{L^2}\Bigr\}_{n\in\mathbb{N}\cup\{0\}},
  \]
  with the eigenvalue $0$ being simple and each eigenvalue
  $\lambda_n = \tfrac{n^2\pi^2}{L^2}$, $n\geq 1$, having multiplicity $2$.
  The continuous spectrum is
  \[
    \sigma_c(\mathcal{H}_Z)
    = [0,+\infty)\setminus\Bigl\{\tfrac{n^2\pi^2}{L^2}\Bigr\}_{n\in\mathbb{N}\cup\{0\}}.
  \]

\item[2)] For $Z < 0$, the spectrum of $\mathcal{H}_Z$ in
  $L^2(\mathcal{G})$ is
  \[
    \sigma(\mathcal{H}_Z)
    = \{\lambda_-\}\cup[0,+\infty),
    \qquad
    \lambda_- = -\frac{N^2}{Z^2},
  \]
  where $\lambda_-$ is a simple eigenvalue with eigenfunction
  $(0,W_{\lambda_-})$, $W_{\lambda_-}=\{e^{\frac{N}{Z}(x-L)}\}_{j=1}^N$,
  and
  \[
    \sigma_p(\mathcal{H}_Z)
    = \{\lambda_-\}
      \cup\Bigl\{\tfrac{n^2\pi^2}{L^2}\Bigr\}_{n\in\mathbb{N}\cup\{0\}},
    \qquad
    \sigma_c(\mathcal{H}_Z)
    = [0,+\infty)\setminus\Bigl\{\tfrac{n^2\pi^2}{L^2}\Bigr\}_{n\in\mathbb{N}\cup\{0\}}.
  \]
\end{enumerate}
\end{proposition}

\begin{proof}
From the boundary conditions in $D(\mathcal{H}_Z)$, the operator  $\mathcal{H}_Z$ decomposes as the orthogonal direct sum
\begin{equation}\label{eq:direct-sum}
  \mathcal{H}_Z
  = \mathcal{H}^{\mathrm{per}} \oplus \mathcal{H}_Z^{\delta'},
\end{equation}
where $\mathcal{H}^{\mathrm{per}}$ is $-d^2/dx^2$
on $L^2([-L,L])$ with domain
\[
  D(\mathcal{H}^{\mathrm{per}})
  = \bigl\{f\in H^2(-L,L)
    \;\big|\; f(-L)=f(L),\;f'(-L)=f'(L)\bigr\},
\]
and $\mathcal{H}_Z^{\delta'}$ is the self-adjoint operator
$-d^2/dx^2$ on $L^2(\dot{\mathcal{G}})$ with domain
\begin{equation}\label{eq:delta-star-star}
  D(\mathcal{H}_Z^{\delta'})
  = \Bigl\{\mathbf{g}\in H^2\bigl(\dot{\mathcal{G}}\bigr)
    \;\Big|\;
    g_1'(L)=\cdots=g_N'(L),\quad
    \sum_{j=1}^N g_j(L)=Z\,g_1'(L)
  \Bigr\},
\end{equation}
where $\dot{\mathcal{G}}$ represents the star graph obtained by removing the
loop from $\mathcal{G}$.
Consequently,
\begin{equation}\label{eq:spectrum-sum}
  \sigma(\mathcal{H}_Z)
  = \sigma(\mathcal{H}^{\mathrm{per}})
    \cup\sigma(\mathcal{H}_Z^{\delta'}).
\end{equation}

\medskip
Since $[-L,L]$ is compact, $\mathcal{H}^{\mathrm{per}}$ has compact
resolvent and hence purely discrete spectrum.
The eigenvalue problem
\begin{equation}\label{eq:per-ev}
  -f'' = \lambda f \quad\text{on }(-L,L),
  \qquad
  f(-L)=f(L),\quad f'(-L)=f'(L),
\end{equation}
has the following solution structure.
\begin{itemize}
  \item $\lambda < 0$: only the trivial solution; hence no negative
        eigenvalues.
  \item $\lambda = 0$: the constant function $f\equiv c$,
        giving a simple eigenvalue.
  \item $\lambda > 0$: general solution
        $f(x)=A\cos(\sqrt{\lambda}\,x)+B\sin(\sqrt{\lambda}\,x)$.
        Imposing periodicity forces $\sqrt{\lambda}=n\pi/L$ for some
        $n\in\mathbb{N}$, and both $A$ and $B$ are free.
        Hence each eigenvalue
        $\lambda_n=\tfrac{n^2\pi^2}{L^2}$ ($n\geq 1$) has
        \emph{multiplicity two}, with eigenfunctions
        $\cos\!\bigl(\tfrac{n\pi}{L}x\bigr)$ and
        $\sin\!\bigl(\tfrac{n\pi}{L}x\bigr)$.
\end{itemize}
Therefore $\sigma(\mathcal{H}^{\mathrm{per}})
=\{n^2\pi^2/L^2\}_{n\geq 0}$, with $\lambda_0=0$ simple and
$\lambda_n$ of multiplicity $2$ for $n\geq 1$.

\medskip
By extension theory applied in $\dot{\mathcal{G}}$
(see Section~6.2 in \cite{AC}):
\begin{itemize}
  \item For $Z>0$: $\sigma(\mathcal{H}_Z^{\delta'})
        =\sigma_{ac}(\mathcal{H}_Z^{\delta'})=[0,+\infty)$ with
        no eigenvalues.
  \item For $Z<0$: $\sigma_{ac}(\mathcal{H}_Z^{\delta'})=[0,+\infty)$
        and $\sigma_p(\mathcal{H}_Z^{\delta'})=\{\lambda_-\}$,
        $\lambda_-=-N^2/Z^2<0$, with eigenfunction
        $W_{\lambda_-}=\{e^{\frac{N}{Z}(x-L)}\}_{j=1}^N\in
         L^2(\dot{\mathcal{G}})$.
\end{itemize}

\medskip
By the decomposition \eqref{eq:direct-sum}, the eigenvalue problem
$\mathcal{H}_Z\mathbf{U}=\lambda\mathbf{U}$ with
$\mathbf{U}=(f,\mathbf{g})=(f,g_1,\ldots,g_N)$ reads as
\begin{equation}\label{eq:spe-ext}
  \begin{cases}
    -f'' = \lambda f,
      & x\in(-L,L), \\[4pt]
    -g_j'' = \lambda g_j,
      & x\in(L,+\infty),\; j=1,\ldots,N, \\[4pt]
    f(-L)=f(L),\quad f'(-L)=f'(L), \\[4pt]
    g_1'(L)=\cdots=g_N'(L), \quad
    \textstyle\sum_{j=1}^N g_j(L)=Z\,g_1'(L).
  \end{cases}
\end{equation}
Note in \eqref{eq:spe-ext} that the two subsystems are independent and any
eigenvector of $\mathcal{H}_Z$ is either of the form $(f,\mathbf{0})$,
$(0,\mathbf{g})$, or a sum of two such independent solutions at the
same eigenvalue. However, a solution with both components nonzero exists only if $\lambda$ is
simultaneously an eigenvalue of $\mathcal{H}^{\mathrm{per}}$ and of
$\mathcal{H}_Z^{\delta'}$.
The eigenvalues of $\mathcal{H}^{\mathrm{per}}$ are non-negative,
and $\mathcal{H}_Z^{\delta'}$ has no eigenvalues in $[0,+\infty)$,
so no such $\lambda$ exists.

\medskip
Finally, let $\lambda<0$ (with $\lambda\neq\lambda_-$ for $Z<0$).
We show $\lambda\in\rho(\mathcal{H}_Z)$.
By the direct sum \eqref{eq:direct-sum},
$(\mathcal{H}_Z-\lambda)^{-1}
 =(\mathcal{H}^{\mathrm{per}}-\lambda)^{-1}
  \oplus(\mathcal{H}_Z^{\delta'}-\lambda)^{-1}$.
Since $\lambda<0$ lies in the resolvent of $\mathcal{H}^{\mathrm{per}}$
(which has non-negative spectrum) and in the resolvent of
$\mathcal{H}_Z^{\delta'}$ (by the known spectral data recalled above),
both components are bounded operators.
Hence $(\mathcal{H}_Z-\lambda)^{-1}$ is bounded, i.e.,
$\lambda\in\rho(\mathcal{H}_Z)$.

\noindent This completes the proof.
\end{proof}

We describe next the resolvent of $\mathcal{H}_Z$.

\begin{proposition}\label{espec1-ext}
  Let $Z\in\mathbb{R}\setminus\{0\}$, $N\geq 1$, and
  $\lambda\in\rho(\mathcal{H}_Z)$.
  The resolvent operator
  $R(\lambda;\mathcal{H}_Z)
   =(\lambda I_{N+1}-\mathcal{H}_Z)^{-1}
   \colon L^2(\mathcal{G})\to D(\mathcal{H}_Z)$
  is given by
  \begin{equation}\label{eq:res-block}
    R(\lambda;\mathcal{H}_Z)
    = \operatorname{diag}\!\bigl(R(\lambda;\mathcal{H}^{\mathrm{per}}),\;
                                R(\lambda;\mathcal{H}_Z^{\delta'})\bigr).
  \end{equation}
  Moreover, for given $k^2=-\lambda>0$ and
  $(f,\{g_j\})\in L^2(\mathcal{G})$, if
  \begin{equation}\label{eq:res-action}
    (\phi,\psi_1,\dots,\psi_N)
    = R(\lambda;\mathcal{H}_Z)(f,g_1,\dots,g_N),
  \end{equation}
  then $\phi(x)=Ae^{kx}+Be^{-kx}+\phi_p(x)$ and
  $\psi_j(x)=A_je^{-kx}+\psi_{j,p}(x)$, where
  \begin{equation}\label{eq:particular}
    \phi_p(x)=\frac{-1}{2k}\!\left[
      \int_{-L}^x e^{k(x-s)}f(s)\,ds
      -\int_x^L e^{-k(x-s)}f(s)\,ds\right],
    \qquad
    \psi_{j,p}(x)=\frac{1}{2k}\int_x^\infty e^{-k(s-x)}g_j(s)\,ds.
  \end{equation}
  The coefficients $[A,B,A_1,\dots,A_N]$ are determined by the
  block-diagonal linear system $\mathbf{M}^{-1}\mathbf{b}$, where
  the matrix $\mathbf{M}$ and right-hand side $\mathbf{b}$ decouple
  as
  \begin{equation}\label{eq:block-M}
    \mathbf{M}
    = \begin{pmatrix}\mathbf{M}^{\mathrm{per}} & 0 \\ 0 & \mathbf{M}^{\delta'}\end{pmatrix},
    \qquad
    \mathbf{b}
    = \begin{pmatrix}\mathbf{b}^{\mathrm{per}}\\\mathbf{b}^{\delta'}\end{pmatrix},
  \end{equation}
  with the loop block $\mathbf{M}^{\mathrm{per}}\in M_2(\mathbb{R})$
  encoding the periodic conditions on $[A,B]$:
  \begin{equation}\label{eq:loop-block}
    \begin{split}
      M_1^{\mathrm{per}} &= \bigl[\sinh(kL),\;-\sinh(kL)\bigr],
      \qquad
      b_1^{\mathrm{per}} = \phi_p(-L)-\phi_p(L),\\
      M_2^{\mathrm{per}} &= \bigl[k\sinh(kL),\;k\sinh(kL)\bigr],
      \qquad\;\;
      b_2^{\mathrm{per}} = \phi_p'(-L)-\phi_p'(L),
    \end{split}
  \end{equation}
  and the half-line block $\mathbf{M}^{\delta'}\in M_N(\mathbb{R})$
  encoding the $\delta'$-conditions on $[A_1,\dots,A_N]$: for
  $j=2,\dots,N$,
  \begin{equation}\label{eq:hl-equal-deriv}
    M_{j-1}^{\delta'}
    = \bigl[\underbrace{ke^{-kL}}_{1\text{st}},0,\dots,0,
            \underbrace{-ke^{-kL}}_{j\text{th}},0,\dots,0\bigr],
    \qquad
    b_{j-1}^{\delta'} = \psi_{j,p}'(L)-\psi_{1,p}'(L),
  \end{equation}
  and for the sum condition,
  \begin{equation}\label{eq:hl-sum}
    M_N^{\delta'}
    = \bigl[e^{-kL}(1+kZ),\;e^{-kL},\;\dots,\;e^{-kL}\bigr],
    \qquad
    b_N^{\delta'}
    = Z\psi_{1,p}'(L)-\sum_{j=1}^N\psi_{j,p}(L).
  \end{equation}
\end{proposition}

\begin{proof}
Let $(f,\mathbf{g})\in L^2(\mathcal{G})$ and suppose
$R(\lambda;\mathcal{H}_Z)(f,\mathbf{g})=(\phi,\boldsymbol{\psi})$.
The block decomposition \eqref{eq:res-block} follows directly from
Proposition~\ref{pointesp-ext} and the direct-sum structure
$\mathcal{H}_Z=\mathcal{H}^{\mathrm{per}}\oplus\mathcal{H}_Z^{\delta'}$
established in \eqref{eq:direct-sum}: the equation
$(\lambda I_{N+1}-\mathcal{H}_Z)(\phi,\boldsymbol{\psi})=(f,\mathbf{g})$
decouples as
\[
  (\lambda-\mathcal{H}^{\mathrm{per}})\phi = f,
  \qquad
  (\lambda I_N-\mathcal{H}_Z^{\delta'})\boldsymbol{\psi} = \mathbf{g},
\]
two independent problems on the loop and on the star graph,
respectively.

For the second statement, one applies the variation-of-parameters
method to each decoupled system. The particular solution $\phi_p$ in \eqref{eq:particular} is the
standard variation-of-parameters formula for
$-\phi''+k^2\phi=f$ on $(-L,L)$.
Imposing the periodic conditions $\phi(-L)=\phi(L)$ and
$\phi'(-L)=\phi'(L)$ on the general homogeneous solution
$Ae^{kx}+Be^{-kx}$ yields the $2\times 2$ system
\eqref{eq:loop-block}.
Row $M_1^{\mathrm{per}}$ arises from $\phi(-L)=\phi(L)$ and
involves only $A-B$; Row $M_2^{\mathrm{per}}$ arises from
$\phi'(-L)=\phi'(L)$ and involves only $A+B$.
The two rows are linearly independent for $k>0$ (since
$\sinh(kL)>0$), so $[A,B]$ is uniquely determined.

\noindent The particular solution $\psi_{j,p}$ in \eqref{eq:particular} is the
variation-of-parameters formula for $-\psi_j''+k^2\psi_j=g_j$ on
$(L,+\infty)$ that decays at infinity.
Imposing the equal-derivative condition $g_j'(L)=g_1'(L)$ for
$j=2,\dots,N$ and the Robin-sum condition $\sum_j g_j(L)=Zg_1'(L)$
on the general solution $A_je^{-kx}+\psi_{j,p}(x)$ yields the
$N\times N$ system \eqref{eq:hl-equal-deriv}--\eqref{eq:hl-sum}.
This system is invertible for all $\lambda\in\rho(\mathcal{H}_Z^{\delta'})$
(see Section~6.2 in \cite{AC}), so $[A_1,\dots,A_N]$ is uniquely
determined.

Since the two blocks share no unknowns, the full $(N+2)\times(N+2)$
matrix $\mathbf{M}$ is block-diagonal as stated in \eqref{eq:block-M},
and its invertibility is immediate.
This finishes the proof.
\end{proof}

\begin{theorem}\label{group-ext}
  Let $Z\in\mathbb{R}\setminus\{0\}$, $N\geq 1$.
  The unitary group generated by $\mathcal{H}_Z$ on $D(\mathcal{H}_Z)$
  acts as
  \begin{equation}\label{eq:group-decomp}
    e^{it\mathcal{H}_Z}(f,\mathbf{g})
    = \bigl(W_1(t)f,\;W_2(t)\mathbf{g}\bigr),
  \end{equation}
  where $W_1(t)f=e^{it\mathcal{H}^{\mathrm{per}}}f$ is the classical
  unitary group acting on $f\in L^2_{\mathrm{per}}([-L,L])$, and
  $W_2(t)\mathbf{g}=e^{-it\mathcal{H}_Z^{\delta'}}\mathbf{g}$ is the
  unitary group generated by $\mathcal{H}_Z^{\delta'}$ on the star
  graph $\dot{\mathcal{G}}$.
\end{theorem}

\begin{proof}
The proof follows from Proposition~\ref{espec1-ext} and the
identification of the resolvent as the Laplace transform of the
group (see \cite[Section~1.7]{pazy}).
Using the block-diagonal resolvent \eqref{eq:res-block} we get
\begin{align*}
  e^{it\mathcal{H}_Z}(f,\mathbf{g})
  &= \frac{1}{2\pi i}
     \int_{\gamma-i\infty}^{\gamma+i\infty}
     e^{\lambda t}\,R(\lambda;\mathcal{H}_Z)(f,\mathbf{g})\,d\lambda \\
  &= \frac{1}{2\pi i}\left(
       \int_{\gamma-i\infty}^{\gamma+i\infty}
       e^{\lambda t}\,R(\lambda;\mathcal{H}^{\mathrm{per}})f\,d\lambda
       \;,\;
       \int_{\gamma-i\infty}^{\gamma+i\infty}
       e^{\lambda t}\,R(\lambda;\mathcal{H}_Z^{\delta'})\mathbf{g}\,d\lambda
     \right) \\
  &= \bigl(W_1(t)f,\;W_2(t)\mathbf{g}\bigr),
\end{align*}
where $\gamma>\beta+1$ with $\beta=0$ for $Z\geq 0$ and
$\beta=N^2/Z^2$ for $Z<0$ (the negative eigenvalue of
$\mathcal{H}_Z^{\delta'}$; see Proposition~\ref{pointesp-ext}).
The integral converges uniformly in $t$ on every compact interval
$[\delta,1/\delta]$, $\delta>0$.
\end{proof}

\begin{theorem}[Local and global well-posedness]\label{global-ext}
  Let $\mathbf{U}_0\in H^1(\mathcal{G})$ and $p>1$.
  Then the integral equation associated with \eqref{NLS},
  \begin{equation}\label{eq:inteq-ext}
    \mathbf{U}(t)
    = e^{it\mathcal{H}_Z}\mathbf{U}_0
      -i\int_0^t e^{i(t-s)\mathcal{H}_Z}
       |\mathbf{U}(s)|^{p-1}\mathbf{U}(s)\,ds,
  \end{equation}
  has a unique solution
  $\mathbf{U}\in C((-T,T);H^1(\mathcal{G}))$ for some
  $T=T(\|\mathbf{U}_0\|_{H^1})>0$ with $\mathbf{U}(0)=\mathbf{U}_0$.
  For any $T_0<T$, the data-solution map
  $\mathbf{U}_0\mapsto\mathbf{U}\in C([-T_0,T_0];H^1(\mathcal{G}))$
  is continuous; for $p>2$ it is of class $C^2$.
  For $1<p<5$ the solution extends globally in time.
\end{theorem}

\begin{proof}
By Theorem~\ref{group-ext}, the integral equation \eqref{eq:inteq-ext}
decouples according to the direct-sum structure
$\mathcal{H}_Z=\mathcal{H}^{\mathrm{per}}\oplus\mathcal{H}_Z^{\delta'}$.
Specifically, the mapping
$J_{\mathbf{U}_0}\colon C([-T,T];H^1(\mathcal{G}))\to
 C([-T,T];H^1(\mathcal{G}))$ given by
\[
  J_{\mathbf{U}_0}[\mathbf{U}](t)
  = e^{it\mathcal{H}_Z}\mathbf{U}_0
    +\int_0^t e^{i(t-s)\mathcal{H}_Z}F(\mathbf{U}(s))\,ds,
  \qquad F(\mathbf{U})=|\mathbf{U}|^{p-1}\mathbf{U},
\]
decomposes component-wise:
\begin{equation}\label{eq:J-decomp}
  J_{\mathbf{U}_0}[\mathbf{U}](t)
  = \bigl(J_{f_0}^{\mathrm{per}}[\phi](t),\;
          J_{\mathbf{g}_0}^{\delta'}[\boldsymbol{\psi}](t)\bigr),
\end{equation}
where $J^{\mathrm{per}}$ is the fixed-point map for the periodic NLS
on the loop driven by $W_1(t)$, and $J^{\delta'}$ is the fixed-point
map for the NLS on the star graph driven by $W_2(t)$.

\medskip
Since $\mathcal{H}^{\mathrm{per}}$ is a self-adjoint operator with
compact resolvent on $L^2_{\mathrm{per}}([-L,L])$, Stone's theorem
yields a unitary group $W_1(t)=e^{it\mathcal{H}^{\mathrm{per}}}$
that preserves $H^1_{\mathrm{per}}([-L,L])$.
Local (and global, for $1<p<5$) well-posedness of the periodic
Cauchy problem
\[
  i\partial_t\phi + \partial_x^2\phi + |\phi|^{p-1}\phi = 0,
  \quad
  \phi(\cdot,0)=f_0\in H^1_{\mathrm{per}}([-L,L]),
\]
in $C([-T,T];H^1_{\mathrm{per}}([-L,L]))$ follows from a standard
contraction-mapping argument in the energy space, using the
Gagliardo--Nirenberg inequality on $[-L,L]$ and conservation of
charge and energy; see \cite{Caz}.

\medskip
The operator $\mathcal{H}_Z^{\delta'}$ defined on $D(\mathcal{H}_Z^{\delta'})$
(see \eqref{eq:delta-star-star}) coincides with the self-adjoint
$\delta'$-interaction Hamiltonian $\mathcal{H}_Z^{\delta'}$ on $\dot{\mathcal{G}}$ studied in \cite{AngGol17b}, with coupling parameter $\lambda=Z$.
Local and global well-posedness of the associated NLS in
$H^1(\dot{\mathcal{G}})$, together with the $H^1$-bound for the
linear evolution $W_2(t)=e^{-it\mathcal{H}_Z^{\delta'}}$, is
established in \cite[Theorem~3.22]{AngGol17b}.

%
%

Since the no vertex conditions are imposed in the energy space
$H^1(\mathcal{G})$, the one-dimensional Sobolev embedding
$W^{1,2}\hookrightarrow L^\infty$ ensures that $F(\mathbf{U})$ is
well-defined in $H^1(\mathcal{G})$ for all $p>1$.

\medskip
By the above, the decomposed map \eqref{eq:J-decomp} is a contraction
on $C([-T,T];H^1(\mathcal{G}))$ for $T=T(\|\mathbf{U}_0\|_{H^1})>0$.
The Banach fixed-point theorem gives a unique solution
$\mathbf{U}\in C([-T,T];H^1(\mathcal{G}))$;
the continuity of the data-solution map follows from standard
estimates, and the $C^2$ regularity for $p>2$ follows from the
Implicit Function Theorem applied to $\Upsilon(\mathbf{W}_0,\mathbf{W})
=\mathbf{W}-J_{\mathbf{W}_0}[\mathbf{W}]$
(see \cite[Corollary~5.6]{LP}).

\medskip
Conservation of charge and energy,
\begin{equation}\label{eq:conserved}
  Q(\mathbf{U})=\|\mathbf{U}(t)\|_{L^2(\mathcal{G})}^2
  = Q(\mathbf{U}_0),
  \qquad
  E_Z(\mathbf{U}(t)) = E_Z(\mathbf{U}_0),
\end{equation}
hold along the flow because $\mathcal{H}_Z$ is self-adjoint on
$D(\mathcal{H}_Z)$.
For $1<p<5$ the Gagliardo--Nirenberg inequality gives control of $\|\mathbf{U}(t)\|_{H^1(\mathcal{G})}$ in terms of conserved
quantities following \cite[Section~2]{AdaNoj14}, so $T$ can be taken
arbitrarily large.
This finishes the proof.
\end{proof}

\section{Existence of coupled dnoidal plus soliton profiles on a looping-edge graph}\label{S6}

In this section, we establish the existence and orbital stability of standing-wave solutions to \eqref{NLS} in the cubic case ($p=3$) on looping-edge graphs with $N$ half-lines where $\mathbf{U}(x,t)=e^{i\omega t}\Theta(x)$, with $\Theta= (\Phi, \Psi_1, \Psi_2, \cdot\cdot\cdot, \Psi_N) \in D'_Z=D'_{Z, N}$, where $D_Z'\subset D(\mathcal{H}_Z)$ (see \eqref{eq:DZ-final}) with
\begin{equation}
    \label{Ndelta'2}
    \begin{split}
        D'_Z := \Big\{ \mathbf{U} \in H^2(\mathcal{G}) \Bigm| & \;\; \phi(L) = \phi(-L), \quad \phi'(L) = \phi'(-L) = \psi_1'(L) = \cdots = \psi_N'(L), \\
        &  \text{and}\;\;\; \sum_{i=1}^N \psi_i(L) = Z \psi_1'(L) \Big\}
    \end{split}
\end{equation}
and $Z \in \mathbb{R}\setminus\{0\}$.

\begin{figure}
    \centering
    \includegraphics[width=0.75\linewidth]{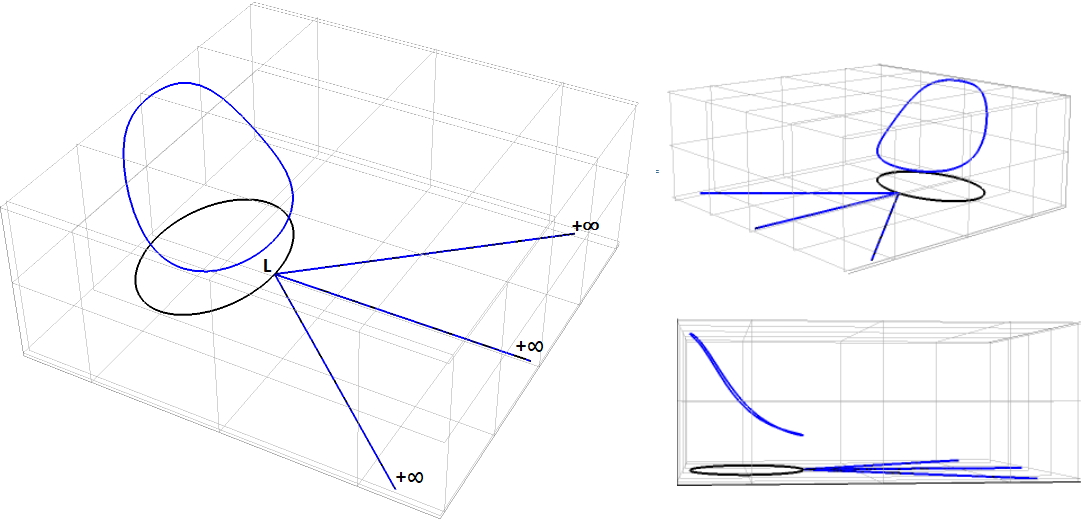}
    \caption{For a looping-edge graph with $N=3$, we exhibit profile elements in $D_Z'$ with $\Phi$ as centered periodic \textit{dnoidal} wave and trivial tail-type profiles on all half-lines.}
    \label{Fig3}
\end{figure}

We note that elements of $D'_Z$ need not be continuous at the vertex $\nu = L$. Here, we are interested in the case where the $\Phi$ is given by the {\it dnoidal} Jacobi elliptic function while each $\Psi_i$ has a tail-type profile (possibly the trivial solution; see Figure~\ref{Fig3}) determined by the classical soliton solution (modulo translation) for the cubic NLS model on the line:
\begin{equation}\label{soliton}
\Psi_i(x) = (2\omega)^{1/2} \, \mathrm{sech}\left(\sqrt{\omega}(x - L) + a_i\right), \quad x \geq L, \quad a_i > 0.
\end{equation}

We point out that equation \eqref{standingN0} is reduced to the system 
   \begin{equation}\label{uv}
 \left\{ \begin{array}{ll}
  -\Phi''(x)+\omega \Phi(x)-|\Phi(x)|^{2} \Phi(x)=0,\;\;\;\;\;\;\;x\in (-L,L),\\
   -\Psi_i''(x)+\omega \Psi_i(x)-|\Psi_i(x)|^{2} \Psi_i(x)=0, \;\;\;x\in (L, +\infty),  \; i=1,...,N\\
(\Phi,\Psi_1, \Psi_2,..., \Psi_N)\in D_{Z}^{\prime}.
  \end{array} \right.
 \end{equation}

Any solution of \eqref{uv} with positive components $(f, \mathbf{g}) \in D'_Z$, where $f$ is even and periodic with period $2L$, must have $\mathbf{g} \equiv 0$. We therefore begin by showing the existence of such periodic profiles.

\begin{proposition}\label{dnoidal}  Let $L>0$ arbitrary but fixed. Then, there exists a smooth mapping of periodic functions $\omega \in (\frac{\pi^ 2}{2L^ 2}, +\infty)\to \Phi_\omega\in H^1_{per}([-L, L])$ such that the periodic profile $\Phi_\omega$ satisfies
 \begin{equation}\label{d0}
 -\Phi_\omega''(x)+\omega \Phi_\omega(x)-\Phi^3_\omega(x)=0,\quad x\in [-L,L],
 \end{equation}
and with a profile of dnoidal type, namely,
  \begin{equation}\label{d1}
\Phi_\omega(x)=\eta_1 dn\Big(\frac{\eta_1}{\sqrt{2}} x;k\Big)
\end{equation}
with the elliptic modulus $k\in (0,1)$ such that
\begin{equation}\label{d2}
k^2(\eta,\omega)=\frac{2\omega-2\eta^2_2}{2\omega-\eta^2_2},\;\;\; \eta^ 2_1+\eta^2_2=2\omega,\;\; 0<\eta_2<\eta_1.
\end{equation}
Moreover, $k=k(\omega)$, $\eta_i=\eta_i(\omega)$ are smooth functions  for $\omega \in (\frac{\pi^ 2}{2L^ 2}, +\infty)$, with the map $\omega\to k(\omega)$ being strictly increasing.
\end{proposition}

\begin{proof} See Theorem 2.1 in \cite{An3}.
\end{proof}

As a scholium to the proof of Proposition~\ref{dnoidal}, it is shown that the function $\varsigma\colon \Omega \to \mathbb{R}$ defined by
\begin{equation}
    \label{fix1}
    \varsigma(\eta, \omega) = \frac{\sqrt{2}}{\sqrt{2\omega - \eta^2}} K(k(\eta, \omega)), \quad \Omega := \left\{ (\eta, \omega) \mid \omega > \frac{\pi^2}{2L^2},\ \eta \in (0, \sqrt{\omega}) \right\}
\end{equation}
satisfies $\partial_\eta \varsigma(\eta, \omega) < 0$, and that for any admissible fixed $\omega$, there exists a unique $\eta(\omega) \in (0, \sqrt{\omega})$ such that $\varsigma(\eta(\omega), \omega) = L$.

As an immediate consequence of Proposition~\ref{dnoidal}, we obtain the following result:

\begin{corollary}\label{d3} 
Let $L>0$ arbitrary but fixed. There is a   smooth mapping $\omega \in (\frac{\pi^ 2}{2L^ 2}, +\infty)\to e^{i\omega t}(\Phi_\omega, \mathbf 0)$ of standing wave solutions for the cubic NLS model on a looping-edge  graph.
\end{corollary}

Let us now study the standing-wave solutions for the cubic \eqref{NLS} where each half-line profile $\Psi_j$ is a non-trivial tail-soliton profile. In this form,  we start our analysis  with the $\Psi_j$'s expression in \eqref{soliton} with a starting point defined by the shift $a_j$.
For each half-line profile the boundary values at the vertex read
\[
\Psi_j(L)=\sqrt{2\omega}\,\operatorname{sech}(a_j),
\qquad
\Psi_j'(L)=-\sqrt{2}\,\omega\,\operatorname{sech}(a_j)\tanh(a_j).
\]
Since $Z<0$ and
\[
\sum_{j=1}^N \Psi_j(L)=Z\Psi_1'(L)>0,
\]
it follows that $\Psi_j'(L)<0$ for all $j$, and therefore $a_j>0$.

The continuity of derivatives at the vertex,
\[
\Psi_1'(L)=\cdots=\Psi_N'(L),
\]
is equivalent to the existence of a common parameter $\mu>0$ such that
\begin{equation}\label{eqmu}
\operatorname{sech}(a_j)\tanh(a_j)=\mu,
\qquad j=1,\dots,N.
\end{equation}
Observe that the function $t\mapsto t\sqrt{1-t^2}$ on $(0,1)$ attains its maximum value $1/2$ at $t=1/\sqrt{2}$, and that for every $\mu\in(0,1/2)$ the equation
\[
t\sqrt{1-t^2}=\mu
\]
admits exactly two solutions,
\[
t_-(\mu)\in\Big(0,\tfrac{1}{\sqrt{2}}\Big)
\quad\text{and}\quad
t_+(\mu)=\sqrt{1-t_-^2(\mu)}\in\Big(\tfrac{1}{\sqrt{2}},1\Big).
\]
Since $\operatorname{sech}(a)=\sqrt{1-\tanh^2(a)}$, equation \eqref{eqmu} therefore admits, for each $j$, exactly two distinct shift solutions,
\[
a_-(\mu)=\operatorname{arctanh}(t_-(\mu))
\quad\text{and}\quad
a_+(\mu)=\operatorname{arctanh}(t_+(\mu)).
\]
Hence, for fixed $\mu$, each half-line shift $a_j$ must belong to the set
$\{a_-(\mu),a_+(\mu)\}$.

\smallskip
Let $n\in\{0,\dots,N\}$ denote the number of half-lines selecting one of the two shifts, and $N-n$ the number selecting the complementary one. Introducing the parameter
\[
\alpha=\alpha(\omega):=-Z\sqrt{\omega}>0,
\]
the sum condition
\begin{equation}\label{sumcond}
\sum_{j=1}^N \Psi_j(L)=Z\Psi_1'(L)
\end{equation}
reduces to the compatibility relation
\[
n\,s_- + (N-n)\,s_+ = \alpha\, s_+ t_+,
\]
where $s_\pm=\operatorname{sech}(a_\pm)$. Since $t_\pm=s_\mp$, this identity can be rewritten as
\begin{equation}\label{eq:kn}
n\,\sqrt{1-t_-^2} + (N-n)\,t_- = \alpha\, t_- \sqrt{1-t_-^2}.
\end{equation}

In particular, when all half-line profiles are continuous at the vertex ($n=0$, $n=N$ or the tadpole case $N=1$), equation \eqref{eq:kn} reduces to an invertible relation $\alpha=\alpha(t_-)$ provided $\alpha>N$, \emph{i.e.}, provided the existence threshold
\[
\omega>\frac{N^2}{Z^2}
\]
is satisfied.

\medskip
For fixed $n$, equation \eqref{eq:kn} admits a solution $t_-\in(0,1)$ if and only if $\alpha$ belongs to the image of the function
\[
\alpha_n(t):=\frac{n}{t}+\frac{N-n}{\sqrt{1-t^2}},
\qquad t\in(0,1),
\]
that is, if and only if
\[
\alpha\ge \alpha_{\min}(n),
\qquad
\alpha_{\min}(n):=\min_{t\in(0,1)} \alpha_n(t).
\]
A direct minimization shows that the minimum is attained at
\[
t^2=\frac{n^{2/3}}{n^{2/3}+(N-n)^{2/3}},
\]
and yields the explicit threshold
\[
\alpha_{\min}(n)=\big(n^{2/3}+(N-n)^{2/3}\big)^{3/2}.
\]
Consequently, half-line configurations with \emph{non-equal shifts} exist if and only if
\[
-Z\sqrt{\omega}\ge \beta(N),
\qquad
\beta(N):=\min_{1\le n\le N-1} \alpha_{\min}(n).
\]
The minimum is attained by the least balanced split of the half-lines, yielding
\begin{equation}\label{betadef}
    \beta(N)=(1+(N-1)^{\tfrac{2}{3}})^{\tfrac{3}{2}}
\end{equation}
In all cases, one must have
\[
\omega\ge \frac{\beta^2(N)}{Z^2}.
\]

In this work, we study the existence and stability of standing waves for the case  in which the tail components satisfy $a_1=a_2=\cdots=a_N\equiv\gamma^*$ and therefore $\Psi_1= \Psi_2=\cdots=\Psi_N\equiv \Psi$ (see Remark~\ref{notequal}, item $(3)$, for some comments about the different shifts case. Moreover, the existence and stability of these profiles will be the focus of a future study). It is not difficult to see that the tail-profile $\Psi=\Psi_{\omega, N,Z}$ is given by the expression
\begin{equation}\label{Ntails}
\Psi_{\omega, Z, N}(x)=\sqrt{2\omega} \, \mathrm{sech}\Big( \sqrt{\omega}(x-L)+\tanh^{-1}\Big(\frac{-N}{Z\sqrt{\omega}}\Big)\Big),\quad x\geq L
\end{equation}
with $\omega>\frac{N^2}{Z^2}$ and $Z<0$. We note that $N\Psi_{\omega, Z, N}(L)=Z \Psi'_{\omega, Z, N}(L)$ (see \eqref{Ndelta'2}).

\begin{figure}
    \centering
    \includegraphics[width=0.5\linewidth]{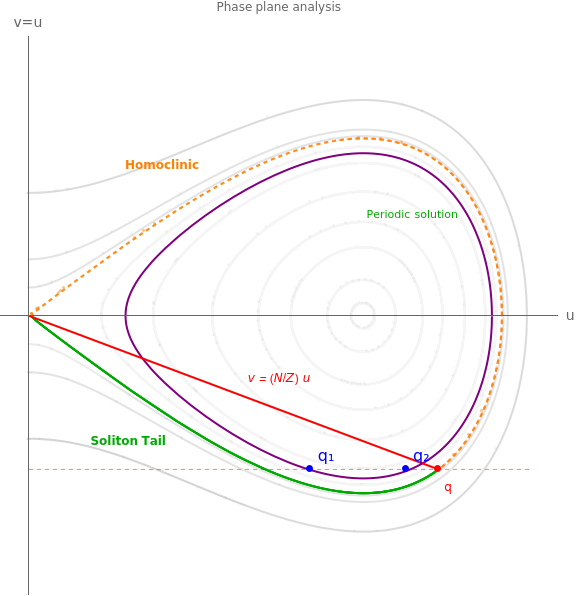}
    \caption{Orange: homoclinic orbit for the standing NLS equation for positive solutions. Green: Soliton tail to be considered. Red: the line $v=\frac{N}{Z}u$. Blue: intersection points of the periodic orbit at the derivative level of $q$. Purple: inner periodic orbit starting and ending at the same $q_i$. Purple + Green: a feasible profile $\Theta\in D'_Z$.  }
    \label{phase}
\end{figure}

\begin{remark}\label{Rem4}
 The condition $\omega>\frac{N^2}{Z^2}$ in the definition of $\Psi_{\omega,Z,N}$ is sharp in the sense that one cannot expect the existence of positive solutions decaying at infinity that solve the second equation in \eqref{uv} when $\omega\le\frac{N^2}{Z^2}$. This can be geometrically understood as the necessity of nontrivial intersections between the homoclinic orbit associated with the stationary cubic NLS equation and the initial value conditions represented as a function of $\Psi(L)$ by the line $\Psi'(L)=\frac{N}{Z}\Psi(L)$. See Figure \ref{phase}.
\end{remark}

Taking into account \eqref{Ntails}, system \eqref{uv} can be translated into
\begin{equation}\label{simp}
 \left\{ \begin{array}{ll}
  -\Phi''(x)+\omega \Phi(x)-\Phi^3(x) =0,\;\;\;\;\;\;\;x\in (-L,L),\\
  \Phi(-L)=\Phi(L), \ \ \Phi'(-L)=\Phi'(L)=\Psi_{\omega,Z,N}'(L)\ \mbox{ and }\ N\Psi_{\omega,Z,N}(L)=Z\Psi_{\omega,Z,N}'(L).
  \end{array} \right.
 \end{equation}
It is worth mentioning that $\Psi(L)$ and $\Psi'(L)$ are no longer unknowns in system \eqref{simp}. In fact,
$$
\Psi(L)=\frac{\sqrt{2}}{-Z}\sqrt{\omega Z^2 - N^2} \quad \text{and} \quad \Psi'(L)=-\frac{N\sqrt{2}}{Z^2}\sqrt{\omega Z^2 - N^2}
$$
See $q=(\Psi(L),\Psi'(L))$ in Figure \ref{phase}. The inflection point of the periodic function must be at or below $v = \Psi'(L)$ in order to guarantee existence of compatible initial conditions in $D'_Z$. The intersection point(s) $q_i$ correspond to suitable initial points of the shifted periodic solution.

Regarding the periodic profile $\Phi$, solution of \eqref{simp}, it is clear that the initial value $\Phi(-L) = \Phi(L)$ is constrained to be such that $\Phi'(L) = \Psi'(L)$. The existence of $\Phi$ depends on the size of $\omega$ relative to $N$ and $Z$. One expects existence only if the shift $\gamma_*$ does not place the initial value too close to an inflection point of $\Psi$. This can be seen geometrically as the necessity of intersection between the inner periodic orbits and the derivative level fixed by the $v$-coordinate of $q$ (see Figure \ref{phase} and \eqref{Interval} below).

Before demonstrating the existence of profiles in $D'_{Z}$ for $Z<0$, we first describe the behavior of an inflection point of the dnoidal wave. Denote $D(x) = \mathrm{dn}(x; k)$, $x \in [-K(k), K(k)]$, with $k \in (0,1)$. We have that $D''(x) = 0$ if and only if $\mathrm{sn}^2(x; k) = \mathrm{cn}^2(x; k)$. Then, from the identity $\mathrm{sn}^2(x; k) + \mathrm{cn}^2(x; k) = 1$ for all $k$, we get $\mathrm{sn}^2(x; k) = \frac{1}{2}$. Thus, from the identity $k^2 \mathrm{sn}^2(x; k) + \mathrm{dn}^2(x; k) = 1$ for all $k$, we arrive at the condition $\mathrm{dn}^2(x; k) = 1 - \frac{k^2}{2}$ for $x \in [-K(k), K(k)]$. 

In this way, there are exactly two symmetric inflection points in $(-K(k), K(k))$, denoted $x_0 = x_0(k) \in (0, K(k)/2)$ and $-x_0$, satisfying
\begin{equation}\label{x0}
\mathrm{dn}^2(x_0; k) = 1 - \frac{k^2}{2}, \quad \text{for every fixed} \;\; k.
\end{equation}
Therefore,
\begin{equation}\label{allinfle}
x_0,\; 2K - x_0,\; 2K + x_0,\; 4K - x_0,\; 4K + x_0,\; \ldots,\;\; -x_0,\; -2K + x_0,\; -2K - x_0,\; -4K + x_0,\; -4K - x_0,\; \ldots
\end{equation}
are all the inflection points of the dnoidal profile on $\mathbb{R}$.

Moreover, it can be seen that $x_0(k) = F\left(\frac{\pi}{4}; k\right)$, with $x_0(k)$ an increasing function of $k \in (0,1)$, satisfying $x_0(1) = F\left(\frac{\pi}{4}; 1\right) = \log(1 + \sqrt{2})$ (see Appendix \ref{A2}).

Recall from Proposition \ref{dnoidal} that the map $\omega \mapsto k(\omega)$ is strictly increasing. For $Z < 0$, $N \ge 1$, denote $k_0 := k\left(\frac{N^2}{Z^2}\right)$ as the infimum value of $k(\omega)$ in the interval $\left(\frac{N^2}{Z^2}, \infty\right)$. Let $r_1 > r_2 > 0$ be the roots of the polynomial $p(x) := \frac{Z^2 k_0^4}{16 N^2}x^2 - x + \frac{N^2}{Z^2}$. Using the quadratic formula, it is not hard to see that $\frac{N^2}{Z^2} < r_2 < \frac{2N^2}{Z^2} < r_1$. Define
\begin{equation}
\label{Interval}
I = I(N, Z, L) := \left( \frac{\pi}{2L^2}, \infty \right) \cap \left( \left( \frac{N^2}{Z^2}, r_2 \right) \cup (r_1, \infty) \right).
\end{equation}

\begin{theorem}\label{existence1}
Let $L>0$, $Z<0$ and $N\geq 1$. There is a smooth mapping $\omega\in I(N,Z,L)\mapsto e^{i\omega t} \Theta_\omega$ of standing wave 
solutions for the cubic NLS model on a looping edge graph in which $\Theta_\omega=(\Phi_{\omega, a},  (\Psi_{\omega, Z, N})_{i=1}^N)\in D'_Z$, where $\Phi_{\omega, a}(x)=\Phi_{\omega}(x-a)$, $x\in [-L,L]$, $a=a(\omega)\in(0,L)$, $\Phi_{\omega}$ is the dnoidal profile in \eqref{d1}, and $\Psi_{\omega, Z, N}$ is given in \eqref{Ntails}.
\end{theorem}

\begin{proof}
Recall $k(\eta,\omega)$ defined in \eqref{d2}. Consider the set  \begin{equation}
    \label{set}
      \mathcal{S}:=\left\{ (\omega,\eta,a)\in \R^3 \ \mid\  \omega\in I, \ \  \eta\in(0,\sqrt{\omega}) \ \ \mbox{and}\ \  L-\frac{\sqrt{2}x_0(k(\eta,\omega))}{\sqrt{2w-\eta^2}}<a<L\right\}
\end{equation}

Define on $\mathcal{S}$ the function $h:\mathcal{S}\to \R^2$ as $h(\omega,\eta,a):=(\varsigma(\eta,\omega),g(\omega,\eta,a))$ where $\varsigma$ is defined on \eqref{fix1} and $g:\mathcal{S}\to \R$ is defined by 
\begin{equation}
g(\omega,\eta,a)=\frac{k(\eta,\omega)^2(2\omega-\eta^2)}{\sqrt{2}}f(\omega,\eta,a)-\sqrt{\omega Z^2-N^2}\frac{N\sqrt{2}}{Z^2},
\end{equation}
where $$f(\omega,\eta,a):=cn\left(\frac{\sqrt{2\omega-\eta^2}}{\sqrt{2}}(L-a);k(\eta,\omega)\right)sn\left(\frac{\sqrt{2\omega-\eta^2}}{\sqrt{2}}(L-a);k(\eta,\omega)\right).$$

Since $dn(u,k)>0$ always and $cn(u,k)<sn(u,k)$ when $u=\frac{\sqrt{2\omega-\eta^2}}{\sqrt{2}}(L-a)\in \left(x_0(k(\eta,\omega),L)\right)$, it is not hard to see that for $\eta_*,\omega_*$ fixed it follows $\partial_a f<0$, which implies $f(\omega_*,\eta_*,a)$ is a strictly decreasing function of $a$. Moreover, $g(\omega_*,\eta_*,a)$ is also strictly decreasing.

Now, note that for any $\omega_*, \eta_*$ fixed, $f(\omega_*,\eta_*,a)\to 0$ as $a\to L$, which implies $g(\omega_*,\eta_*,a)$ tends to a strictly negative quantity as $a\to L$. Similarly, $f(\omega_*,\eta_*,a)\to 1/2$ as $a\to L-\frac{\sqrt{2}}{\sqrt{2\omega-\eta^2}}x_0(k(\eta,\omega))$, which implies that $g(\omega_*,\eta_*,a)$ tends to 
\begin{equation}\label{620}
\frac{k(\eta,\omega)^2(2\omega-\eta^2)}{2\sqrt{2}}-\sqrt{\omega Z^2-N^2}\frac{N\sqrt{2}}{Z^2}
\end{equation}
as $a\to L-\frac{\sqrt{2}}{\sqrt{2\omega-\eta^2}}x_0(k(\eta,\omega))$.

Let us see that \eqref{620} is a positive quantity. In fact,
\begin{equation}
    \begin{split}
    \frac{k(\eta,\omega)^2(2\omega-\eta^2)}{2\sqrt{2}}>\sqrt{\omega Z^2-N^2}\frac{N\sqrt{2}}{Z^2} \longleftrightarrow \frac{Z^2k(\eta,\omega)^4(2\omega-\eta^2)^2}{16N^2}+\frac{N^2}{Z^2}>\omega.
    \end{split}
\end{equation}

Since $k_0\le k(\eta,\omega)$ and $\eta<\sqrt{\omega}$ it is enough to have
\begin{equation}\label{ri}
\frac{Z^2k(\eta,\omega)^4(2\omega-\eta^2)^2}{16N^2}+\frac{N^2}{Z^2}>\frac{Z^2k_0^4w^2}{16N^2}+\frac{N^2}{Z^2}>w.
\end{equation}
The rightmost inequality on \eqref{ri} is equivalent to have $p(\omega)>0$ which correspond to the choice $\omega\in (\frac{N^2}{Z^2},r_2)\cup(r_1,\infty)$ in the definition of $I$.

From the Intermediate value theorem and since $g(\omega_*,\eta_*,a)$ is strictly decreasing, there exists a unique $a=a(\eta_*,\omega_*)$ such that $g_(\omega_*,\eta_*,a)=0$.

Let $\omega_0\in I$ be arbitrary. From the scholium of Proposition \ref{dnoidal}, there exists a unique $\eta_0\in (0,\sqrt{\omega_0})$ such that $\varsigma(\eta_0,\omega_0)=2L$. Furthermore, denote with $a_0$ the unique value such that $(\omega_0,\eta_0,a_0)\in \mathcal{S}$ and $g(\omega_0,\eta_0,a_0)=0$. Note in particular $h(\omega_0,\eta_0,a_0)=(2L,0)$. Since both, $\partial_\eta \varsigma(\eta_0,\omega_0)$ and $\partial_a g(\omega_0,\eta_0,a_0)$ are different than $0$, the matrix $$J_{h,(\eta,a)}(\omega_0,\eta_0,a_0)=\begin{pmatrix}
    \partial_\eta \varsigma(\omega_0,\eta_0,a_0) & 0 \\ \partial_\eta g(\omega_0,\eta_0,a_0) & \partial_a g(\omega_0,\eta_0,a_0)
\end{pmatrix}$$
is invertible regardless the value of $\partial_\eta g(\omega_0,\eta_0,a_0)$. By the Implicit function theorem, there is a neighborhood $\mathcal{O}(\omega_0)\subset I$ and a unique smooth function $\Lambda:\mathcal{O}(\omega_0)\to \R^2$ such that $\Lambda(\omega_0)=(\eta_0,a_0)$ and for any $\omega\in\mathcal{O}(\omega_0)$ it follows $h(\omega,\Lambda(\omega))=(2L,0)$. 
Besides, since $\omega_0\in I$ is arbitrary and $\Lambda$ is unique, $\mathcal{O}(\omega_0)$ can be taken to be $I$.

The associated dnoidal wave solution $\Phi_{\omega,a(\omega)}$ with $\eta_1\equiv\eta_1(\omega), \eta_2\equiv\eta_2(\omega)$ has fundamental period $2L$ and the mapping $\omega\in I\to \Phi_{\omega,a}\in H^1_{\mbox{per}}([-L,L])$ is smooth.

Finally, $\Theta_\omega\in D'_Z$ given that $N\Psi_{\omega,Z,N}(L)=Z\Psi'_{\omega,Z,N}(L)$ and since the fact $g(\omega,\eta(\omega),a(\omega))=0$ is equivalent to $\Psi'_{\omega,Z,N}(L)=\Phi_{\omega,a(\omega)}'(L)$.
\end{proof}

\begin{remark}\label{notequal} From the proof of Theorem \ref{existence1} it is worth mentioning that
\begin{enumerate}
    \item It is clear from Figure \ref{phase} the existence of periodic solutions with orbit tangent to the line $v=\Psi'(L)$, \textit{i.e}, periodic solutions starting from its inflection point; however the argument won't guarantee smoothness since the Jacobian matrix in the last two variables $J_{h,(\eta,a)}$ is not invertible when $\frac{\sqrt{2\omega-\eta^2}}{\sqrt{2}}(L-a)$ is at the inflection point $x_0$. 
    
    \item The proof can be easily adapted to obtain $\frac{\sqrt{2\omega-\eta^2}}{\sqrt{2}}(L-a)\in (0,x_0)$ (\textit{i.e.,} to pick $q_2$ instead of $q_1$ in Figure \ref{phase}) without significative modifications.

    \item By exploiting the existence of tail profiles \eqref{soliton} with distinct shifts \(a_j\) that satisfy the \(\delta'\)-coupling conditions on the star graph \(\dot{\mathcal{G}}\) introduced at the beginning of this subsection, we may follow a strategy analogous to that used in the proof of Theorem~\ref{existence1} to construct a smooth curve \(\omega \mapsto \Theta_\omega \in D'_Z\) of standing-wave solutions for the cubic NLS. The stability properties of these profiles will be investigated in a future work.

\end{enumerate}
\end{remark}

The values of $\omega, \ N$ and $Z$ modify strongly the behavior of the profile $\Psi=\Psi_{\omega, Z, N}$ in \eqref{Ntails}. For instance, as $\frac{N}{Z}$ gets close to $0$, the value of $\Psi'(L)$ gets closer to $0$ implying that the shift $a$ is close to $L$ and that most of the feasible periodic solutions satisfy that its maximum value is below the maximum value $\Psi(L)$ of $\Psi$. On the other hand, as $|Z|$ goes to $N\sqrt{\frac{2}{\omega}}$ by the right, the point $q$ becomes closer to the inflection point of $\Psi$ implying that the shift $a$ is close to $x_0$ and most of the feasible periodic solutions satisfy that its maximum value is above the maximum value $\Psi(L)$ of $\Psi$ (see Figure \ref{Fig4}).

\section{(In)stability analysis}

Our first goal in this section is to study the stability properties of the standing-wave
solutions established in Corollary~\ref{d3}.
To this end, we first identify the symmetry group of the NLS model
\eqref{NLS} on a looping-edge graph and then provide the appropriate
definition of orbital stability.

We begin by observing that the NLS model \eqref{NLS} with $(-\Delta,D(\mathcal{H}_Z))$ and posed on a
looping-edge graph admits \emph{two} continuous symmetries.

\begin{itemize}
  \item Phase invariance: If $\mathbf{U}=(\phi,\psi_1,\ldots,\psi_N)$
        is a solution of \eqref{NLS}, then so is
        $e^{i\theta}\mathbf{U}$ for any $\theta\in\mathbb{R}$,
        since the equation and the boundary conditions in
        $D(\mathcal{H}_Z)$ are invariant under global phase rotation.
  \item Loop translation invariance: The periodic
        boundary conditions $\phi(-L)=\phi(L)$ and
        $\phi'(-L)=\phi'(L)$ on the loop component are invariant
        under spatial translation by any $s\in\mathbb{R}/(2L\mathbb{Z})$.
        Consequently, if $(\phi(x),\psi_1(x),\ldots,\psi_N(x))$
        is a solution of \eqref{NLS}, then so is
        $(\phi(x+s),\psi_1(x),\ldots,\psi_N(x))$ for any such $s$.
        Note that no analogous translation is available on the
        half-lines, whose profiles are pinned at the vertex $\nu=L$
        by the vertex conditions in $D(\mathcal{H}_Z)$.
\end{itemize}

These two symmetries generate the group
$G := \mathbb{R}\times\mathbb{R}/(2L\mathbb{Z})$
acting on $H^1(\mathcal{G})$ by
\begin{equation}\label{eq:group-action}
  T_{(\theta,s)}\mathbf{U}
  := e^{i\theta}\bigl(\phi(\cdot+s),\,\psi_1,\ldots,\psi_N\bigr),
  \qquad (\theta,s)\in G.
\end{equation}
The \emph{orbit} of a profile $\Theta=(\Phi,\Psi_1,\ldots,\Psi_N)$
under $G$ is
\begin{equation}\label{eq:orbit}
  \mathcal{O}(\Theta)
  := \bigl\{T_{(\theta,s)}\Theta
     :\theta\in\mathbb{R},\;s\in\mathbb{R}/(2L\mathbb{Z})\bigr\}
   = \bigl\{e^{i\theta}(\Phi(\cdot+s),\Psi_1,\ldots,\Psi_N)
     :\theta\in\mathbb{R},\;s\in\mathbb{R}/(2L\mathbb{Z})\bigr\}.
\end{equation}
The symmetries \eqref{eq:group-action} naturally motivate the
following definition of orbital stability (see \cite{GSS2, GSS1}).

\begin{definition}\label{dsta}
  The standing wave $\mathbf{U}(x,t)=e^{i\omega t}\Theta(x)$ is said
  to be \textit{orbitally stable} in a Banach space $X$ if for any
  $\varepsilon>0$ there exists $\eta>0$ with the following property:
  if $\mathbf{U}_0\in X$ satisfies
  $\|\mathbf{U}_0-\Theta\|_X<\eta$, then the solution $\mathbf{U}(t)$
  of \eqref{NLS} with $\mathbf{U}(0)=\mathbf{U}_0$ exists for all
  $t\in\mathbb{R}$ and
  \[
    \sup_{t}\,
    \inf_{(\theta, s)}\,
    \bigl\|\mathbf{U}(t)
           - e^{i\theta}\bigl(\Phi(\cdot+s),\Psi_1,\ldots,\Psi_N\bigr)
    \bigr\|_X
    < \varepsilon.
  \]
  Equivalently, $\sup_{t\in\mathbb{R}}
  \mathrm{dist}_X(\mathbf{U}(t),\mathcal{O}(\Theta))<\varepsilon$,
  where $\mathcal{O}(\Theta)$ is the orbit \eqref{eq:orbit}.
  Otherwise, the standing wave is said to be
  \textit{orbitally unstable} in $X$.
\end{definition}

The space $X$ in Definition \ref{dsta} for the model \eqref{NLS} will be the energy-space $H^ 1(\mathcal{G})$, endowed with the action $-\Delta$ defined on $D(\mathcal{H}_Z)$.

Now, for a fixed $\omega > 0$, let $\mathbf{U}(x,t) = e^{i\omega t}\Theta(x)$, with $\Theta = (\Phi, \Psi_1, \Psi_2, \ldots, \Psi_N)$ a stationary solution of \eqref{NLS} such that $\Theta \in D(\mathcal{H}_Z)$. Then, for the action functional
\begin{equation}\label{S1}
\mathbf{S}(\mathbf{U}) = E_Z(\mathbf{U}) + \omega Q(\mathbf{U}), \qquad \mathbf{U} \in H^1(\mathcal{G}),
\end{equation}
we have $\mathbf{S}'(\Theta) = \mathbf{0}$.

Next, for $\mathbf{U} = \mathbf{U}_1 + i \mathbf{U}_2$ and $\mathbf{W} = \mathbf{W}_1 + i \mathbf{W}_2$, where the vector functions $\mathbf{U}_j$ and $\mathbf{W}_j$ ($j = 1, 2$) are assumed to have real components, it is not hard to verify that the second variation of $\mathbf{S}$ at $\Theta$ is given by
\begin{equation}\label{S2}
\mathbf{S}''(\Theta)(\mathbf{U}, \mathbf{W}) 
= \langle \mathcal{L}_{+,Z} \mathbf{U}_1, \mathbf{W}_1 \rangle 
+ \langle \mathcal{L}_{-,Z} \mathbf{U}_2, \mathbf{W}_2 \rangle,
\end{equation}
where the two $(N+1)\times (N+1)$ diagonal operators $\mathcal{L}_{\pm,Z}$ are defined by
\begin{equation}\label{L+}
\begin{aligned}
\mathcal{L}_{+,Z} &= \mathrm{diag}\Big(-\partial_x^2 + \omega - p|\Phi|^{p-1},\,
-\partial_x^2 + \omega - p|\Psi_1|^{p-1},\, \ldots,\,
-\partial_x^2 + \omega - p|\Psi_N|^{p-1} \Big),\\[0.5em]
\mathcal{L}_{-,Z} &= \mathrm{diag}\Big(-\partial_x^2 + \omega - |\Phi|^{p-1},\,
-\partial_x^2 + \omega - |\Psi_1|^{p-1},\, \ldots,\,
-\partial_x^2 + \omega - |\Psi_N|^{p-1} \Big).
\end{aligned}
\end{equation}

We note that the operators $\mathcal{L}_{\pm,Z}$ defined above are self-adjoint with domain $D(\mathcal{L}_{\pm,Z}) \equiv D(\mathcal{H}_Z)$. Moreover, since $\Theta \in D'_Z$ satisfies the system \eqref{uv}, we have $\mathcal{L}_{-,Z}\Theta = \mathbf{0}$, and therefore the kernel of $\mathcal{L}_{-,Z}$ is nontrivial.

\subsection{Stability of dnoidal+trivial soliton profiles on looping edge graphs}

In this subsection, we establish the orbital stability of each member of the family of standing wave solutions $e^{i\omega t}(\Phi_\omega, \mathbf{0})$ described in Corollary~\ref{d3}. In this context, it is known from \cite{GSS2, GSS1} that the Morse index and the nullity of the operators $\mathcal{L}_{\pm,Z}$ in \eqref{L+}, evaluated around the profile $(\Phi_\omega, \mathbf{0})$, are fundamental in determining the orbital stability of standing wave solutions for NLS models (see Theorem~\ref{2main} in the Appendix).

%
%
%
%
%

\begin{theorem}\label{Mtad-ext}
  Let $(\mathcal{L}_{+,Z}, D(\mathcal{H}_Z))$ be the self-adjoint
  operator in \eqref{L+} for fixed $Z\in\mathbb{R}\setminus\{0\}$,
  determined by the dnoidal profile $(\Phi_\omega,\mathbf{0})$ of
  Proposition~\ref{dnoidal} for $\omega>\frac{\pi^2}{2L^2}$, namely
  \begin{equation}\label{spectrumL+-ext}
    \mathcal{L}_{+,Z}
    = \operatorname{diag}\!\bigl(
        -\partial_x^2+\omega-3\Phi_\omega^2,\,
        -\partial_x^2+\omega,\,\ldots,\,
        -\partial_x^2+\omega
      \bigr).
  \end{equation}
  Then:
  \begin{enumerate}
    \item[1)] The kernel of $\mathcal{L}_{+,Z}$ satisfies
      \begin{itemize}
        \item For $Z>0$ or $Z<0$ with $\omega> N^2/Z^2$:
              $\ker(\mathcal{L}_{+,Z})$ is one-dimensional and
              spanned by $(\Phi'_\omega,\mathbf{0})$.
        \item For $Z<0$ and $\omega=N^2/Z^2$: $\ker(\mathcal{L}_{+,Z})$
              is two-dimensional and spanned by
              $(\Phi'_\omega,\mathbf{0})$ and $(0,W_{\lambda_-})$,
              where $W_{\lambda_-}=\{e^{\frac{N}{Z}(x-L)}\}_{j=1}^N$.
      \end{itemize}
    \item[2)] In case $Z>0$ or $Z<0$ with $\omega> N^2/Z^2$, the Morse index of $\mathcal{L}_{+,Z}$ is exactly one.
  \end{enumerate}
\end{theorem}

\begin{proof}
  We recall the spectral properties of
  $\mathcal{L}_{dn}=-\partial_x^2+\omega-3\Phi_\omega^2$ on
  $H^2_{\mathrm{per}}([-L,L])$. We have
  $\ker(\mathcal{L}_{dn})=\operatorname{span}\{\Phi'_\omega\}$, $n(\mathcal{L}_{dn})=1$,
  and the unique negative eigenvalue $\lambda_0$ with even eigenfunction
  $\chi_0$ satisfying $\chi_0'(L)=0$.

 \medskip
  In what comes to item $1)$, assume $(f,\mathbf{g})\in\ker(\mathcal{L}_{+,Z})$ with
  $(f,\mathbf{g})\in D(\mathcal{H}_Z)$. Then
$\mathcal{L}_{dn}f=0$ gives $f=\alpha\Phi'_\omega$
  for some $\alpha\in\mathbb{R}$.
  Note $D(\mathcal{H}_Z)$ imposes only the periodic conditions
  $f(-L)=f(L)$ and $f'(-L)=f'(L)$ on the loop component, both satisfied by $\Phi'_\omega$.

On the other hand, $(-\partial_x^2+\omega)g_j=0$ on $[L,+\infty)$
  yields $g_j(x)=A_je^{-\sqrt{\omega}(x-L)}$ (the unique $L^2$
  solution).
  The $D(\mathcal{H}_Z)$ vertex conditions give
  \begin{itemize}
    \item Equal derivatives:
          $-\sqrt{\omega}A_j=-\sqrt{\omega}A_1$ for all $j$,
          hence $A_j=A_1=:A$.
    \item Sum condition:
          $NA=Z(-\sqrt{\omega}A)$, i.e.\ $A(N+Z\sqrt{\omega})=0$.
  \end{itemize}
  Since the loop and half-line conditions in $D(\mathcal{H}_Z)$ are
  completely independent, $\alpha$ and $A$ can be determined separately, namely
  \begin{itemize}
    \item If $Z>0$ or $\omega\neq N^2/Z^2$ (i.e.\
          $N+Z\sqrt{\omega}\neq 0$): then $A=0$, so
          $\mathbf{g}\equiv\mathbf{0}$, and $\alpha$ is free.
          Hence $\ker(\mathcal{L}_{+,Z})=\operatorname{span}\{(\Phi'_\omega,\mathbf{0})\}$.
    \item If $Z<0$ and $\omega=N^2/Z^2$ (i.e.\
          $N+Z\sqrt{\omega}=0$): then $A$ is free as well.
          Hence $\ker(\mathcal{L}_{+,Z})
          =\operatorname{span}\{(\Phi'_\omega,\mathbf{0}),(0,W_{\lambda_-})\}$.
  \end{itemize}

  \medskip
Let us now focus on the Morse index. Since
  $\langle\mathcal{L}_{+,Z}(\Phi_\omega,\mathbf{0}),
   (\Phi_\omega,\mathbf{0})\rangle
  =-2\int_{-L}^L\Phi_\omega^3\,dx<0$,
  the min-max principle gives $n(\mathcal{L}_{+,Z})\geq 1$. Suppose $(f,\mathbf{g})\in D(\mathcal{H}_Z)$
  satisfies $\mathcal{L}_{+,Z}(f,\mathbf{g})=\lambda(f,\mathbf{g})$
  for some $\lambda<0$. From the loop component we have $\mathcal{L}_{dn}f=\lambda f$ with $\lambda<0$. The only option is $\lambda=\lambda_0$ and
  $f=\beta\chi_0$.
  Since $\chi_0$ is even, $\chi_0'(L)=0$.

On the other hand, $(-\partial_x^2+\omega)g_j=\lambda g_j$ on
  $[L,+\infty)$, so $g_j(x)=B_je^{-\sqrt{\omega-\lambda}(x-L)}$.
  The $D(\mathcal{H}_Z)$ conditions give $B_j=B_1=:B$ and
  \begin{equation}\label{eq:B-condition}
    B\bigl(N+Z\sqrt{\omega-\lambda}\bigr) = 0.
  \end{equation}
  Since $\lambda<0$, we have $\omega-\lambda>\omega>0$. If $Z>0$ then $N+Z\sqrt{\omega-\lambda}>0$ and if $Z<0$ and $\omega>N^2/Z^2$ we have $\sqrt{\omega-\lambda}>\sqrt{\omega}>N/|Z|$, hence
  $Z\sqrt{\omega-\lambda}<-N$, giving
  $N+Z\sqrt{\omega-\lambda}<0$.
  In either case $N+Z\sqrt{\omega-\lambda}\neq 0$, so \eqref{eq:B-condition}
  forces $B=0$, hence $\mathbf{g}\equiv\mathbf{0}$.

  Therefore every negative eigenvector has the form
  $(f,\mathbf{g})=(\beta\chi_0,\mathbf{0})$, and the corresponding
  eigenvalue is $\lambda_0$.
\end{proof}


\begin{corollary}\label{dnoidalext-ext}
  Let $a\in(0,L)$ and consider the translated dnoidal profile
  $\Phi_{\omega,a}(x)=\Phi_\omega(x-a)$, $x\in[-L,L]$, for
  $\omega>\frac{\pi^2}{2L^2}$.
  The operator
  \[
    \mathcal{L}_{dn,a}=-\partial_x^2+\omega-3\Phi_{\omega,a}^2
  \]
  on $H^2_{\mathrm{per}}([-L,L])$ has the same spectrum as
  $\mathcal{L}_{dn}$; in particular
  $\ker(\mathcal{L}_{dn,a})=\operatorname{span}\{\Phi'_{\omega,a}\}$ and
  $n(\mathcal{L}_{dn,a})=1$, with negative the same eigenvalue $\lambda_0$ and eigenfunction
  $\zeta_{\lambda_0}(x)=\chi_0(x-a)$.
\end{corollary}

\begin{figure}
    \centering
    \includegraphics[width=\linewidth]{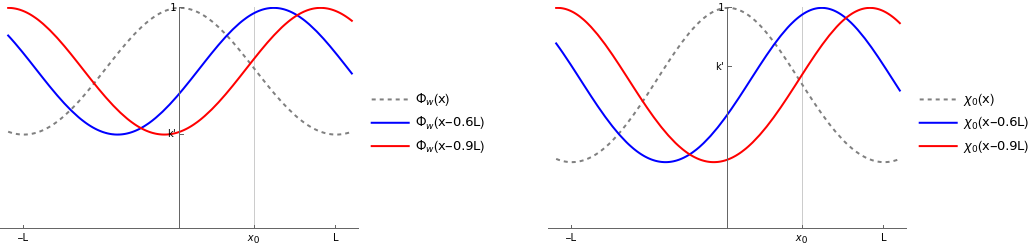}
    \caption{Plots of $\chi_0$ and $\Phi_w$ with sample translations with $a=0.6L$ and $a=0.9L$ using $k=\frac{1}{2}$.}
    \label{fig:transl}
\end{figure}

\begin{figure}
    \centering
    \includegraphics[width=0.75\linewidth]{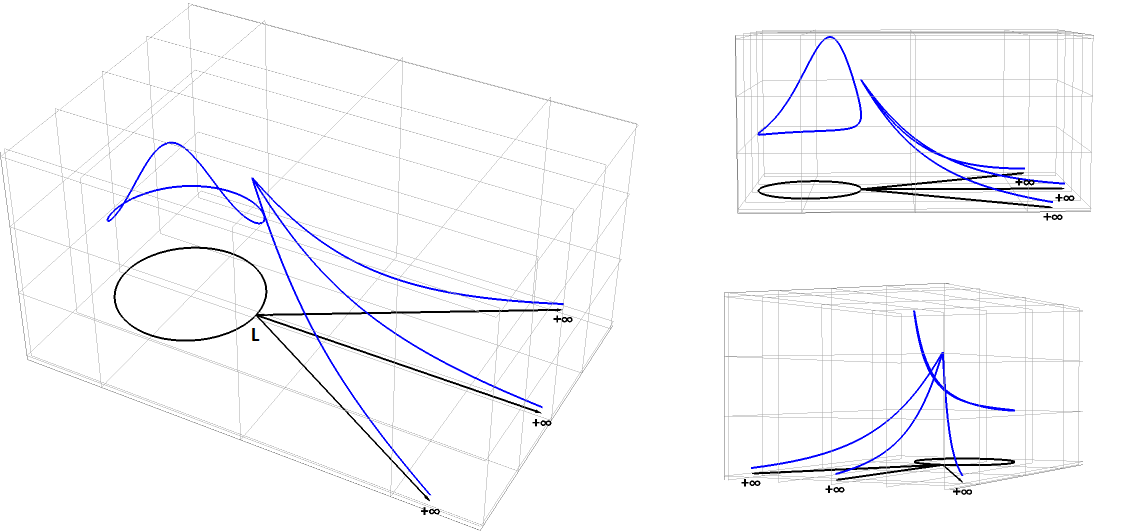}
    \caption{For $N=3$, a sketch of a profile in $D_Z'$ with $Z<0$ where $\Phi_{\omega,a}$ is a shifted periodic \textit{dnoidal} wave and $\psi_1=\psi_2=\psi_3$ are tail-type profiles.}
    \label{Fig4}
\end{figure}

\begin{proof}
    The spectral equivalence $\sigma(\mathcal{L}_{dn,a})=\sigma(\mathcal{L}_{dn})$
  follows from the translation isometry $T_a(f)(x)=f(x+a)$ on
  $H^2_{\mathrm{per}}([-L,L])$. Note the proof of Theorem~\ref{Mtad-ext} applies verbatim with $\Phi_\omega$ replaced
  by $\Phi_{\omega,a}$ and $\chi_0$ replaced by $\chi_0(\cdot-a)$.
  In particular, for all $a\in(0,L)$:
  \begin{equation}\label{eq:zeta-nonzero}
    \zeta'_{\lambda_0}(L) = \chi_0'(L-a) \neq 0,
  \end{equation}
  since the zero of $\chi_0'$ on $[0,L]$ occurs only at the
  endpoints $0$ and $L$ (see Figure~\ref{fig:transl}).
\end{proof}

%
%
%
%

\begin{theorem}\label{Mtad2-ext}
  Let $(\mathcal{L}_{-,Z},D(\mathcal{H}_Z))$ be the self-adjoint
  operator in \eqref{L+} for $Z\in\mathbb{R}\setminus\{0\}$,
  associated with the dnoidal profile $(\Phi_\omega,\mathbf{0})$ of
  Proposition~\ref{dnoidal} for $\omega>\frac{\pi^2}{2L^2}$, namely
  \begin{equation}\label{spectrumL--ext}
    \mathcal{L}_{-,Z}
    = \operatorname{diag}\!\bigl(
        -\partial_x^2+\omega-\Phi_\omega^2,\;
        -\partial_x^2+\omega,\;\ldots,\;
        -\partial_x^2+\omega
      \bigr).
  \end{equation}
  Then the following hold:
  \begin{enumerate}
    \item[1)] The kernel of $\mathcal{L}_{-,Z}$ satisfies
      \begin{itemize}
        \item For $Z>0$ or $Z<0$ with $\omega\neq N^2/Z^2$:
              $\ker(\mathcal{L}_{-,Z})=\operatorname{span}\{(\Phi_\omega,\mathbf{0})\}$.
        \item[2)] For $Z<0$ and $\omega=N^2/Z^2$:
              $\ker(\mathcal{L}_{-,Z})=\operatorname{span}\{(\Phi_\omega,\mathbf{0}),(0,W_{\lambda_-})\}$,
              where $W_{\lambda_-}=\{e^{\frac{N}{Z}(x-L)}\}_{j=1}^N$.
      \end{itemize}
    \item[(2)] For $Z>0$ or $Z<0$ with $\omega> N^2/Z^2$ the operator $\mathcal{L}_{-,Z}$ is non-negative.
  \end{enumerate}
\end{theorem}

\begin{proof}
Consider the operator $\mathcal{L}_0=-\partial_x^2+\omega-\Phi_\omega^2$
on $H^2_{\mathrm{per}}([-L,L])$.
Since $\mathcal{L}_0\Phi_\omega=0$ and $\Phi_\omega>0$, Floquet theory
gives $\ker(\mathcal{L}_0)=\operatorname{span}\{\Phi_\omega\}$ and $\mathcal{L}_0\geq 0$.

\medskip
Let $(f,\mathbf{g})\in\ker(\mathcal{L}_{-,Z})$ with
$(f,\mathbf{g})\in D(\mathcal{H}_Z)$. From the loop component, $\mathcal{L}_0 f=0$ gives $f=\alpha\Phi_\omega$ for some
$\alpha\in\mathbb{R}$. On the half-lines, $(-\partial_x^2+\omega)g_j=0$ on $[L,+\infty)$
gives $g_j(x)=A_je^{-\sqrt{\omega}(x-L)}$.
The $D(\mathcal{H}_Z)$ vertex conditions again yield $A=A_j$ for all $j$ and $A(N+Z\sqrt{\omega})=0$.

For $Z>0$ or $Z<0$ with $\omega\neq N^2/Z^2$, we have
$N+Z\sqrt{\omega}\neq 0$, hence $A=0$ and $\mathbf{g}\equiv\mathbf{0}$.
For $Z<0$ with $\omega=N^2/Z^2$, $A$ is free and
$\mathbf{g}=A\,W_{\lambda_-}$.
The kernel conclusions follow.

In what comes to the non-negativity, the loop part satisfies $\mathcal{L}_0\geq 0$ by Floquet theory.
The half-line part satisfies $-\partial_x^2+\omega\geq 0$ on each
$[L,+\infty)$ for $\omega>\tfrac{N^2}{Z^2}$.
Since $D(\mathcal{H}_Z)$ imposes no coupling between the two components,
we have
\[
  \langle\mathcal{L}_{-,Z}(f,\mathbf{g}),(f,\mathbf{g})\rangle
  = \langle\mathcal{L}_0 f,f\rangle_{L^2([-L,L])}
    + \sum_{j=1}^N\langle(-\partial_x^2+\omega)g_j,g_j\rangle_{L^2([L,+\infty))}
  \geq 0
\]
for all $(f,\mathbf{g})\in D(\mathcal{H}_Z)$, with equality if and only
if $(f,\mathbf{g})\in\ker(\mathcal{L}_{-,Z})$.
\end{proof}

We are in position to prove our stability result.

\begin{theorem}\label{1stability1-ext}
  Let $L>0$ be arbitrary but fixed and $Z\in\mathbb{R}\setminus\{0\}$.
  For every $\omega>\frac{\pi^2}{2L^2}$ (with $\omega$ also satisfying $\omega>N^2/Z^2$ if
  $Z<0$), the standing-wave solution $e^{i\omega t}(\Phi_\omega,\mathbf{0})$
  is orbitally stable in $H^1(\mathcal{G})$ under the flow of the cubic
  NLS equation on a looping-edge graph generated by $\mathcal{H}_Z$ on
  $D(\mathcal{H}_Z)$.
  Here $\Phi_\omega$ is the dnoidal profile \eqref{d1}--\eqref{d2}.
\end{theorem}

\begin{proof}
By Corollary~\ref{d3}, Theorem~\ref{global-ext},
Theorems~\ref{Mtad-ext} and~\ref{Mtad2-ext}, it remains to verify the
non-degeneracy slope condition of the abstract orbital stability
theorem (Theorem~\ref{2main}) applied with the symmetry group
$G=U(1)\times\mathbb{R}/(2L\mathbb{Z})$ acting by
\eqref{eq:group-action}.

\medskip
The symmetry group $G=U(1)\times\mathbb{R}/(2L\mathbb{Z})$ has Lie
algebra $\mathfrak{g}\cong\mathbb{R}^2$. For $(\alpha,\beta)\in\mathfrak{g}$, the augmented action functional and
its associated scalar function are
\begin{equation}\label{eq:d-two-param}
  S_{(\alpha,\beta)}(u)
  = E(u)
    + \frac{\alpha}{2}\|u\|^2_{L^2(\mathcal{G})}
    - \frac{\beta}{2}\operatorname{Im}\!\int_{-L}^{L}\bar{f}\,\partial_x f\,dx,
  \qquad
  d(\alpha,\beta) = S_{(\alpha,\beta)}(\phi_{(\alpha,\beta)}),
\end{equation}
where $\phi_{(\alpha,\beta)}=(e^{i\frac{\beta}{2}x}\Psi_\mu,\mathbf{0})$
is the critical point of $S_{(\alpha,\beta)}$ and
$\mu=\alpha-\frac{\beta^2}{4}$.
Since we work with the standing-wave family
$e^{i\omega t}(\Phi_\omega,\mathbf{0})$, we set $\beta=0$ throughout,
so that $\phi_{(\omega,0)}=(\Phi_\omega,\mathbf{0})$ and write
$d(\omega):=d(\omega,0)$.
The gradient of $d$ encodes two conserved quantities,
\[
  \nabla d(\alpha,\beta)\big|_{\beta=0}
  =\bigl(M(\Phi_\omega),\,-P(\Phi_\omega)\bigr),
\]
where $M(\Phi_\omega)=\frac{1}{2}\|\Phi_\omega\|^2_{L^2([-L,L])}$
is the mass and
\[
  P(\Phi_\omega)
  =\frac{1}{2}\operatorname{Im}\!\int_{-L}^{L}
   \bar\Phi_\omega\,\partial_x\Phi_\omega\,dx = 0
\]
is the momentum, which vanishes since $\Phi_\omega$ is real.
Writing $m(\omega):=M(\Phi_\omega)$ for brevity, from $\partial_\alpha M = m'(\omega)$ and $\partial_\beta^2 d =
-\frac{1}{2}M(\Phi_\omega)=-\frac{m(\omega)}{2}$, and noting that the
off-diagonal entries $\partial_\alpha\partial_\beta
d=-\frac{\beta}{2}m'(\omega)$ and $\partial_\alpha P=\frac{\beta}{2}m'(\omega)$ vanish at $\beta=0$, we obtain the diagonal matrix
\begin{equation}\label{eq:hessian-diag}
  d''(\omega,0)
  =\begin{pmatrix}
      m'(\omega) & 0\\[4pt]
      0 & -\dfrac{m(\omega)}{2}
   \end{pmatrix}.
\end{equation}
In particular, whenever $m'(\omega)>0$, the matrix
\eqref{eq:hessian-diag} is indefinite with exactly one positive
eigenvalue $m'(\omega)$ and one negative eigenvalue
$-m(\omega)/2$.

\medskip
\noindent Let us now prove that in fact $m'(\omega)>0.$ Since the half-line component is identically zero,
$\|(\Phi_\omega,\mathbf{0})\|^2_{L^2(\mathcal{G})}
 =\|\Phi_\omega\|^2_{L^2([-L,L])}$.
By Proposition~\ref{dnoidal} and the theory of elliptic functions
(see \cite[Theorem~3.5]{An3}),
\begin{equation}\label{eq:L2-norm}
  m(\omega)
  = \frac{1}{2}\int_{-L}^{L}\Phi_\omega^2(x)\,dx
  = \frac{2}{L}\,K(k)\,E(k),
\end{equation}
where $K$ and $E$ denote the complete elliptic integrals of first and
second kind, respectively \cite{ByrdFriedman} and the map $k\in(0,1)\mapsto K(k)E(k)$ is strictly increasing.
Since $\frac{dk}{d\omega}>0$ by
Proposition~\ref{dnoidal}, the map $\omega\mapsto m(\omega)$ is
strictly increasing, so $m'(\omega)>0$ for all admissible $\omega$, confirming that the non-degeneracy slope
condition of Theorem~\ref{2main} is satisfied.
This completes the proof.
\end{proof}

\subsection{Stability of dnoidal+tail soliton profiles on looping edge graphs}

\vspace{3mm}

Let us now consider the (in)stability of the standing waves associated with the profiles  $\Theta_\omega=(\Phi_{\omega, a},  (\Psi_{\omega, Z, N})_{i=1}^N)$ determined in Theorem \ref{existence1}. We present a direct proof of the instability. 

\begin{theorem}
\label{thm:instability-HZ}
  Let $L>0$, $Z<0$, $N\geq 1$, and $\omega > N^2/Z^2$.
  The standing wave $e^{i\omega t}\Theta_\omega$, with profile
  $\Theta_\omega \in D_Z'$, is \emph{orbitally unstable} in
  $H^1(\mathcal{G})$ under the NLS flow generated by
  $(\mathcal{H}_Z, D(\mathcal{H}_Z))$.
\end{theorem}

\begin{proof}
For every $\varepsilon>0$ define
\begin{equation}\label{eq:Ueps}
    \mathbf{U}_\varepsilon(0)
    :=\bigl(\Phi_{\omega,a},\;
\Psi_{\omega+\varepsilon,Z,N},\ldots,\Psi_{\omega+\varepsilon,Z,N}\bigr).
  \end{equation}
  Then $\mathbf{U}_\varepsilon(0)\in D(\mathcal{H}_Z)$ for all $\varepsilon>0$. 
Since every component of $\mathbf{U}_\varepsilon(0)$ is an exact
stationary profile for the cubic NLS at its respective frequency
(the loop at $\omega$, all tails at $\omega+\varepsilon$),
the unique $H^1(\mathcal{G})$ solution with initial datum
$\mathbf{U}_\varepsilon(0)$ is
\begin{equation}\label{eq:exact-sol-HZ}
  \mathbf{U}_\varepsilon(t)
  = \bigl(
      e^{i\omega t}\Phi_{\omega,a},\;
      e^{i(\omega+\varepsilon)t}\Psi_{\omega+\varepsilon,Z,N},\;
      \ldots,\;
      e^{i(\omega+\varepsilon)t}\Psi_{\omega+\varepsilon,Z,N}
    \bigr).
\end{equation}
In particular at the time $t_\varepsilon:=\pi/\varepsilon$ one has
$e^{i(\omega+\varepsilon)t_\varepsilon}
 =e^{i\omega t_\varepsilon}\cdot e^{i\pi}
 =-e^{i\omega t_\varepsilon}$,
so all half-line components flip their sign relative to the loop
component simultaneously
\begin{equation}\label{eq:sign-flip-HZ}
  \mathbf{U}_\varepsilon(t_\varepsilon)
  = e^{i\omega t_\varepsilon}
    \bigl(\Phi_{\omega,a},\;
     -\Psi_{\omega+\varepsilon,Z,N},\;\ldots,\;
     -\Psi_{\omega+\varepsilon,Z,N}\bigr).
\end{equation}

\medskip
Since $\omega\mapsto\Psi_{\omega,Z,N}$ is $C^1$ into
$H^1([L,+\infty))$, the
mean-value inequality gives
\begin{equation}\label{eq:init-dist-HZ}
  \|\mathbf{U}_\varepsilon(0)-\Theta_\omega\|_{H^1(\mathcal{G})}^2
  = N\,\|\Psi_{\omega+\varepsilon,Z,N}-\Psi_{\omega,Z,N}
         \|_{H^1([L,+\infty))}^2
  \leq N\,C_\omega^2\,\varepsilon^2,
\end{equation}
where
$C_\omega := \sup_{\omega'\in(\omega,\omega+1)}
             \|\partial_{\omega'}\Psi_{\omega',Z,N}
             \|_{H^1([L,+\infty))}<\infty$.
Hence, as $\varepsilon\to 0$ one has
$$\mathrm{dist}_{H^1}(\mathbf{U}_\varepsilon(0),\mathcal{O}(\Theta_\omega))
 \leq \sqrt{N}\,C_\omega\,\varepsilon\longrightarrow 0.$$ Here $\mathcal{O}(\Theta_\omega)$ denotes the orbit \eqref{eq:orbit}.

\medskip
Let
$\mathbf{V}=e^{i\theta}(\Phi_{\omega,a}(\cdot+s),
 \Psi_{\omega,Z,N},\ldots,\Psi_{\omega,Z,N})
 \in\mathcal{O}(\Theta_\omega)$ be any orbit element
and set $\mu:=\theta-\omega t_\varepsilon$.
Expanding in $L^2(\mathcal{G})$ gives
\begin{align}
  \|\mathbf{U}_\varepsilon(t_\varepsilon)-\mathbf{V}\|_{L^2(\mathcal{G})}^2
  &=
    \|\Phi_{\omega,a}-e^{i\mu}\Phi_{\omega,a}(\cdot+s)\|_{L^2([-L,L])}^2
    + N\,\|\Psi_{\omega+\varepsilon,Z,N}+e^{i\mu}\Psi_{\omega,Z,N}
          \|_{L^2([L,+\infty))}^2 \notag\\
  &= A_N(\varepsilon)+2\cos(\mu)\bigl(P_N(\varepsilon)-q(s)\bigr),
  \label{eq:quad-form-HZ}
\end{align}
where
\begin{align*}
  A_N(\varepsilon)
  &:= 2\|\Phi_{\omega,a}\|_{L^2}^2
     +N\bigl(\|\Psi_{\omega+\varepsilon,Z,N}\|_{L^2}^2
             +\|\Psi_{\omega,Z,N}\|_{L^2}^2\bigr),\\
  P_N(\varepsilon)
  &:= N\langle\Psi_{\omega+\varepsilon,Z,N},
              \Psi_{\omega,Z,N}\rangle_{L^2([L,+\infty))},\\
  q(s)
  &:= \langle\Phi_{\omega,a},\Phi_{\omega,a}(\cdot+s)
             \rangle_{L^2([-L,L])}.
\end{align*}
Minimising \eqref{eq:quad-form-HZ} over $\mu\in[0,2\pi)$ and then
over $s\in\mathbb{R}/(2L\mathbb{Z})$ yields
\begin{equation}\label{eq:inf-orbit-HZ}
  \mathrm{dist}_{L^2}\bigl(\mathbf{U}_\varepsilon(t_\varepsilon),
    \mathcal{O}(\Theta_\omega)\bigr)^2
  \geq A_N(\varepsilon)-2\sup_{s}|P_N(\varepsilon)-q(s)|.
\end{equation}
From Lemma~\ref{lem:lower-bound-HZ} below, there is a constant
$c_{0,N}>0$, independent of $\varepsilon$, such that the right-hand
side is at least $c_{0,N}^2$.
Since $H^1(\mathcal{G})\hookrightarrow L^2(\mathcal{G})$,
\begin{equation}\label{eq:lower-bound-HZ}
  \mathrm{dist}_{H^1}\bigl(\mathbf{U}_\varepsilon(t_\varepsilon),
    \mathcal{O}(\Theta_\omega)\bigr)\geq c_{0,N}>0.
\end{equation}

\medskip
Finally, given $\delta>0$ arbitrary choose $\varepsilon(\delta)>0$ with
$\sqrt{N}\,C_\omega\varepsilon<\delta$.
Then \eqref{eq:init-dist-HZ} gives the estimate
$\mathrm{dist}_{H^1}(\mathbf{U}_\varepsilon(0),\mathcal{O}(\Theta_\omega))<\delta$,
while \eqref{eq:lower-bound-HZ} gives 
$\mathrm{dist}_{H^1}(\mathbf{U}_\varepsilon(t_\varepsilon),\mathcal{O}(\Theta_\omega))
 \geq c_{0,N}>0$.
This is
precisely the negation of orbital stability.
\end{proof}

\begin{lemma}\label{lem:lower-bound-HZ}
  There exists a constant $c_{0,N}>0$, depending only on
  $\omega,L,k,N,Z$, such that
  \begin{equation}\label{eq:lb-HZ}
    A_N(\varepsilon)-2\sup_{s\in\mathbb{R}/(2L\mathbb{Z})}
    |P_N(\varepsilon)-q(s)|
    \geq c_{0,N}^2
    \qquad\text{for all }\varepsilon>0.
  \end{equation}
\end{lemma}

\begin{proof}
Set the following positive constants, all independent of $\varepsilon$:
\[
  \alpha    := \|\Phi_{\omega,a}\|_{L^2([-L,L])}^2,\qquad
  \beta_0   := \|\Psi_{\omega,Z,N}\|_{L^2([L,+\infty))}^2,\qquad
  q_{\min}  := \inf_{s\in\mathbb{R}/(2L\mathbb{Z})} q(s)>0,
\]
where $q_{\min}>0$ because $\Phi_{\omega,a}>0$.
Note $q_{\max}:=\sup_s q(s)=\alpha$.

\medskip
First note the maps
$\varepsilon\mapsto\|\Psi_{\omega+\varepsilon,Z,N}\|_{L^2}^2$
and
$\varepsilon\mapsto
 \langle\Psi_{\omega+\varepsilon,Z,N},\Psi_{\omega,Z,N}\rangle_{L^2}$
are continuous, both with limits equal to $\beta_0$ as
$\varepsilon\to 0^+$.
Write
$\beta_\varepsilon:=\|\Psi_{\omega+\varepsilon,Z,N}\|_{L^2}^2$
and $p_\varepsilon:=\langle\Psi_{\omega+\varepsilon,Z,N},
\Psi_{\omega,Z,N}\rangle_{L^2}$,
so that
\[
  A_N(\varepsilon)=2\alpha+N(\beta_\varepsilon+\beta_0),\qquad
  P_N(\varepsilon)=Np_\varepsilon.
\]
By the Cauchy--Schwarz inequality,
\begin{equation}\label{eq:CS-tail}
  0 < p_\varepsilon
    \leq \sqrt{\beta_\varepsilon\beta_0}
    \leq \tfrac12(\beta_\varepsilon+\beta_0),
  \qquad\text{so}\quad
  2Np_\varepsilon\leq N(\beta_\varepsilon+\beta_0).
\end{equation}

\medskip
Since $Np_\varepsilon$  is independent on $s$ and $q(s)\in[q_{\min},\alpha]$ for all $s$, the supremum in
\eqref{eq:lb-HZ} is attained at one of the two endpoints
\[
  \sup_s|Np_\varepsilon-q(s)|
  = \max\!\bigl(|Np_\varepsilon-q_{\min}|,\;
               |Np_\varepsilon-\alpha|\bigr).
\]
We argue by cases.

\medskip
\noindent\textit{Case 1: $Np_\varepsilon\leq q_{\min}$.}
The supremum equals $\alpha-Np_\varepsilon$. Then
\[
  A_N(\varepsilon)-2\sup_s|P_N-q|
  = N(\beta_\varepsilon+\beta_0)+2Np_\varepsilon
  \geq N\beta_0.
\]

\noindent\textit{Case 2: $Np_\varepsilon\geq\alpha$.}
The supremum equals $Np_\varepsilon-q_{\min}$.
Using \eqref{eq:CS-tail},
\[
  A_N(\varepsilon)-2\sup_s|P_N-q|
  = 2\alpha+N(\beta_\varepsilon+\beta_0)-2Np_\varepsilon+2q_{\min}
  \geq 2\alpha+2q_{\min}.
\]

\noindent\textit{Case 3: $q_{\min}<Np_\varepsilon<\alpha$.} By considering the subcases $Np_\varepsilon-q_{\min}\geq\alpha-Np_\varepsilon$ or $\alpha-Np_\varepsilon\geq Np_\varepsilon-q_{\min}$ the estimates coincide to the same computations done in cases 1 or 2. 

\medskip
In all cases the lower bound is at least
$\min(N\beta_0,\;2\alpha+2q_{\min})=:c_{0,N}^2>0.$
\end{proof}

\appendix
\section{\hspace{1mm}}
\subsection{Orbital (in)stability criterion}
In the Grillakis-Shatah-Strauss (\cite{GSS2,GSS1}) setting we assume
the existence of $C^2$-conserved functionals $E_Z:H^1(\mathcal{G})\to\R$
(energy) with $E_Z(\mathbf{U}(t))=E_Z(\mathbf{U}(0))$, and
$Q:L^2(\mathcal{G})\to\R$ (mass) with
$Q(\mathbf{U}(t))=Q(\mathbf{U}(0))=\|\mathbf{U}(0)\|^2_2$.
In the present setting, the NLS on $\mathcal{G}$ is invariant under
the action of the two-dimensional Lie group
$G=U(1)\times\mathbb{R}/(2L\mathbb{Z})$, given by simultaneous phase
rotation and loop translation,
\begin{equation}\label{eq:group-action_app}
  T_{(\theta,s)}\mathbf{U}
  =\bigl(e^{i\theta}f(\cdot-s),\,e^{i\theta}g_1,\,\ldots,\,e^{i\theta}g_N\bigr),
  \qquad (\theta,s)\in\mathbb{R}/2\pi\mathbb{Z}\times\mathbb{R}/2L\mathbb{Z}.
\end{equation}
The energy $E_Z$ is invariant under the full action \eqref{eq:group-action_app}.
In the applications of this paper, the relevant bound states take the
form $T_{(\omega t,0)}\Theta_\omega$, i.e.\ pure standing waves with
$\beta=0$; the loop-translation parameter is retained in the abstract
framework for completeness but plays no role in the stability
computation beyond determining the orbit.

\smallskip
Suppose that the Cauchy problem associated to \eqref{NLS} is globally
well-posed in the energy space $H^1(\mathcal{G})$, more precisely, it
is guaranteed existence and uniqueness of solutions with continuous
data--solution map.
We suppose the existence of a $C^1$ map on $\mathcal{O}\subset\R$,
\[
  \mathcal{O}\ni\omega\mapsto\Theta_\omega\in D(\mathcal{H}_Z),
\]
of stationary solutions that are critical points of the augmented
action functional
\[
  S_\omega = E_Z + {\omega}Q,
\]
which coincides with $S_{(\omega,0)}$ in the two-parameter notation
\eqref{eq:d-two-param} upon setting $\beta=0$.

\bigskip
The scalar function associated to the two-parameter family is
$d(\omega,\beta)=S_{(\omega,\beta)}(\phi_{(\omega,\beta)})$.
Since $\beta=0$ throughout, we write $d(\omega):=d(\omega,0)$ and note
that the Hessian $d''(\omega,0)$ is the $2\times 2$ diagonal matrix
\eqref{eq:hessian-diag}, whose single negative eigenvalue corresponds
to the translation direction and whose positive eigenvalue corresponds
to the phase direction. Define on $\mathcal{O}$ the function
\begin{equation}
    \rho(\omega_0):=\begin{cases}
        1,& \mbox{if } \partial_\omega \|\Theta_\omega\|^2 >0
            \ \ \mbox{at }\omega=\omega_0,\\
        0, & \mbox{if } \partial_\omega \|\Theta_\omega\|^2 <0
            \ \ \mbox{at }\omega=\omega_0.
    \end{cases}
\end{equation}
For the (in)stability study of $\Theta_\omega$ the main information
will be given by the second variation
\[
  S''(\Theta_\omega)(\mathbf{U},\mathbf{V})
  =\langle\mathcal{L}_{+,Z}\mathbf{U}_1,\mathbf{V}_1\rangle
   +\langle\mathcal{L}_{-,Z}\mathbf{U}_2,\mathbf{V}_2\rangle,
  \qquad
  \mathbf{U}=\mathbf{U}_1+i\mathbf{U}_2,\quad
  \mathbf{V}=\mathbf{V}_1+i\mathbf{V}_2,
\]
in the following sense:

\begin{theorem}\label{2main}
  Suppose $\ker(\mathcal{L}_{-,Z})=\operatorname{span}\{\Theta_\omega\}$
  and that $\ker(\mathcal{L}_{+,Z})$ is spanned by the loop-translation
  generator $(\Theta'_\omega,\mathbf{0})$. Assume also the Morse indices
  $n(\mathcal{L}_{\pm,Z})$ are finite while the rest of the spectra are
  bounded away from $0$. Let
  $\mathcal{H}:=\operatorname{diag}(\mathcal{L}_{+,Z},\mathcal{L}_{-,Z})$
  and let $p = n(\mathcal{H}) - n(d''(\omega,0))$ account for the
  negative directions absorbed by the two-dimensional orbit.
  Then the following hold:
    \begin{enumerate}
        \item If $n(\mathcal{H})=\rho(\omega)=1$, then the standing
              wave $e^{i\omega t}\Theta_\omega$ is orbitally stable in
              $H^1(\mathcal{G})$, with orbit
              \[
                \mathcal{O}(\Theta_\omega)
                =\bigl\{e^{i\theta}(\Theta_\omega(\cdot+s),\mathbf{0})
                 :\theta\in\mathbb{R},\,
                  s\in\mathbb{R}/(2L\mathbb{Z})\bigr\}.
              \]
        \item If $n(\mathcal{H})-\rho(\omega)$ is odd, then the
              standing wave $e^{i\omega t}\Theta_\omega$ is orbitally
              unstable in $H^1(\mathcal{G})$.
    \end{enumerate}
\end{theorem}

\subsection{Jacobian elliptic functions}\label{A2}

We summarize standard definitions and identities for elliptic integrals and Jacobi elliptic functions. For further details see \cite{ByrdFriedman}.

\subsubsection*{A.2.1 Normal Elliptic Integrals of the First and Second Kind}

The \emph{incomplete elliptic integral of the first kind} is defined as
\[
F(\phi; k) = \int_0^\phi \frac{d\theta}{\sqrt{1 - k^2 \sin^2 \theta}},
\]
and the \emph{complete elliptic integral of the first kind} is $K(k) =F\left(\frac{\pi}{2}; k\right).$

The \emph{incomplete elliptic integral of the second kind} is
\[
E(\phi, k) = \int_0^\phi \sqrt{1 - k^2 \sin^2 \theta} \, d\theta,
\]
with the \emph{complete elliptic integral of the second kind} given by $
E(k) = E\left(\frac{\pi}{2}, k\right).
$ Here \( k \in [0,1) \) is the \emph{modulus}, and \( k' = \sqrt{1 - k^2} \) is the \emph{complementary modulus}.

\subsubsection*{A.2.2 Jacobi Elliptic Functions}

Let \( u = F(\phi, k) \) be the elliptic argument. The \emph{Jacobi elliptic functions} \( \operatorname{sn}(u, k) \), \( \operatorname{cn}(u, k) \), and \( \operatorname{dn}(u, k) \) are defined as:
\[
\operatorname{sn}(u, k) = \sin \phi, \quad
\operatorname{cn}(u, k) = \cos \phi, \quad
\operatorname{dn}(u, k) = \sqrt{1 - k^2 \sin^2 \phi}.
\]

We have used in Section \ref{S6} that $F\left(\frac{\pi}{4};k\right)$ is a solution of $sn^2(u,k)=1/2$. The latter is clear given that $\mbox{sin}(\frac{\pi}{4})=\frac{\sqrt{2}}{2}$. We also used that $F\left(\frac{\pi}{4};1\right)=\log{(1+\sqrt{2})}$, which is explicitly in \cite[111.04]{ByrdFriedman}. We also used the following identities:
\[
\operatorname{sn}^2(u,k) + \operatorname{cn}^2(u,k) = 1, \quad
k^2 \operatorname{sn}^2(u,k) + \operatorname{dn}^2(u,k) = 1.
\]

\subsubsection*{A.2.3 Derivatives of Jacobi Elliptic Functions}

Derivatives with respect to the argument \( u \) are
\[
\frac{d}{du} \operatorname{sn}(u,k) = \operatorname{cn}(u,k) \operatorname{dn}(u,k),
\]
\[
\frac{d}{du} \operatorname{cn}(u,k) = -\operatorname{sn}(u,k) \operatorname{dn}(u,k),
\]
and
\[
\frac{d}{du} \operatorname{dn}(u,k) = -k^2 \operatorname{sn}(u,k) \operatorname{cn}(u,k).
\]

Partial derivatives with respect to the modulus \( k \) include (see \cite[710.53]{ByrdFriedman}):
\[
\frac{\partial}{\partial k} \operatorname{dn}(u, k) = \frac{k}{k'^2} \left[ E(\phi, k) - k'^2 u - \frac{\operatorname{sn}(u,k) \operatorname{cn}(u,k)}{\operatorname{dn}(u,k)} \right],
\]
where \( \phi = \operatorname{am}(u,k) \) is the \emph{Jacobi amplitude}.

\subsubsection*{A.2.4 Derivative of the Complete Elliptic Integral of the First Kind}
We have used the following derivatives involving the complete integrals of first and second kind. 
\begin{equation*}
\begin{split}
    \frac{dK}{dk}=\frac{E-{k'}^2K}{k{k'}^2},\ \  \frac{dE}{dk}=\frac{E-K}{k},\\
    \frac{d^2E}{dk^2}=-\frac{1}{k}\frac{dK}{dk}=-\frac{E-{k'}^2K}{k2{k'}^2}.
    \end{split}
\end{equation*}
See, for instance, \cite[710.00, 710.02]{ByrdFriedman}.

\subsubsection*{A.2.5 Properties of the Dnoidal Function}

The function \( \operatorname{dn}(u, k) \) is even and periodic with real period $T = 2K(k)$, so that
\[
\operatorname{dn}(u + 2K(k), k) = \operatorname{dn}(u, k), \quad \operatorname{dn}(-u, k) = \operatorname{dn}(u, k).
\]

It attains its \emph{maximum} value of 1 at \( u = 0 \), and local minima at \( u = K(k) \), where $
\operatorname{dn}(K(k), k) = k'.
$

The second derivative of \( \operatorname{dn}(u,k) \) is
\[
\frac{d^2}{du^2} \operatorname{dn}(u,k) = -k^2 \left[ \operatorname{cn}^2(u,k) - \operatorname{sn}^2(u,k) \right] \operatorname{dn}(u,k).
\]
Since \( \operatorname{dn}(u,k) > 0 \), concavity is governed by the sign of \( \operatorname{cn}^2 - \operatorname{sn}^2 \). Thus, \( \operatorname{dn}(u,k) \) is
\begin{itemize}
    \item \emph{Concave downward} near \( u = 0 \), where \( \operatorname{cn}(u,k) \approx 1 \),
    \item \emph{Concave upward} near \( u = K(k) \), where \( \operatorname{sn}(u,k) \approx 1 \).
\end{itemize}


\section*{Acknowledgments}
This study is financed, in part, by the S\~ao Paulo Research Foundation (FAPESP), Brazil. Process number 2024/20623-7. 

J. Angulo was partially funded by CNPq/Brazil Grant. 

\vskip0.1in
 \noindent
{\bf Data availability statement}. Data supporting this study are included within the article and/or supporting materials.


\begin{thebibliography}{99}
    
\bibitem{AdaNoj15} 
R. Adami, C. Cacciapuoti, D. Finco and D. Noja. 
{\it Stable standing waves for a NLS on star graphs as local minimizers of the constrained energy,} 
J. Differential Equations, 260, 7397--7415. 2016.

\bibitem{AdaNoj14} 
R. Adami, C. Cacciapuoti, D. Finco and D. Noja. 
{\it Variational properties and orbital stability of standing waves for NLS equation on a star graph,} 
J. Differential Equations, 257, 3738--3777. 2014.

\bibitem{AST} 
R. Adami, E. Serra and P. Tilli. 
{\it NLS ground states on graphs,} 
Calc. Var. Partial Differential Equations, 54(1), 743--761. 2015.

\bibitem{ASTcri} 
R. Adami, E. Serra and P. Tilli. 
{\it Negative energy ground states for the $L^2$-critical NLSE on metric graphs,} 
Comm. Math. Phys., 352, 387--406. 2017.

\bibitem{ASTmul} 
R. Adami, E. Serra and P. Tilli. 
{\it Multiple positive bound states for the subcritical NLS equation on metric graphs,} 
Calc. Var., 58(1), 5, 16pp. 2019.


\bibitem{An1} 
J. Angulo. 
{\it Stability theory for two-lobe states on the tadpole graph for the NLS equation,} 
Nonlinearity, 37, 045015. 2024.

\bibitem{An2} 
J. Angulo. 
{\it Stability theory for the NLS on looping-edge graphs,} 
Math. Z., 308, 19. 2024.


\bibitem{An3} 
J. Angulo. 
{\it Nonlinear stability of periodic traveling wave solutions to the Schr\"odinger and the modified Korteweg-de Vries equations,} 
JDE, 235, 1--30. 2007.

\bibitem{AC} 
J. Angulo and M. Cavalcante. 
{\it Nonlinear Dispersive Equations on Star Graphs,} 
$32^o$ Col\'oquio Brasileiro de Matem\'atica, IMPA. 2019.

\bibitem{AC1} 
J. Angulo and M. Cavalcante. 
{\it Linear instability of stationary solitons for the Korteweg-de Vries equation on a star graph,} 
Nonlinearity, 34, 3373--3410. 2021.


\bibitem{AngGol17a} 
J. Angulo and N. Goloshchapova. 
{\it On the orbital instability of excited states for the NLS equation with the $\delta$-interaction on a star graph,} 
Discrete Contin. Dyn. Syst. A., 38(10), 5039--5066. 2018.

\bibitem{AngGol17b} 
J. Angulo and N. Goloshchapova. 
{\it Extension theory approach in the stability of the standing waves for the NLS equation with point interactions on a star graph,} 
Adv. Differential Equations, 23(11--12), 793--846. 2018.

\bibitem{ALN}  J. Angulo,  O. Lopes, and A. Neves, {\it Instability of travelling waves for weakly coupled KdV systems}, Nonlinear Anal. 69 (2008), no. 5-6, pp. 1870--1887. 

\bibitem{AM1}
J. Angulo and A. Mu\~noz,
{\it Airy and Schr\"odinger-type equations on looping-edge graphs}.
arXiv:2410.11729, 2024.%





\bibitem{ANata1} J. Angulo and F. Natali, {\it On the instability of periodic waves for dispersive equations}, Differ. Integral Equ. 29 (2016), no. 9-10, pp. 837--874.






\bibitem{AP1} 
J. Angulo and R. Plaza. 
{\it Instability of static solutions of the sine-Gordon equation on a $Y$-junction graph with $\delta$-interaction,} 
J. Nonlinear Sci., 31, 50. 2021.

\bibitem{AP2} 
J. Angulo and R. Plaza. 
{\it Instability theory of kink and anti-kink profiles for the sine-Gordon on Josephson tricrystal boundaries,} 
Physica D, 427, 133020. 2021.

\bibitem{AP3} 
J. Angulo and R. Plaza. 
{\it Unstable kink and anti-kink profiles for the sine-Gordon on a $Y$-junction graph with $\delta'$-interaction at the vertex,} 
Math. Z., 300, 2885--2915. 2022.

\bibitem{Ardila} 
A. H. Ardila. 
{\it Orbital stability of standing waves for supercritical NLS with potential on graphs,} 
Appl. Anal., 99(8), 1359--1372. 2020.

\bibitem{BK} 
G. Berkolaiko and P. Kuchment. 
{\it Introduction to Quantum Graphs,} 
Math. Surv. Monogr., 186, Amer. Math. Soc., Providence, RI. 2013.

\bibitem{BMP} 
G. Berkolaiko, J. L. Marzuola and D. E. Pelinovsky. 
{\it Edge-localized states on quantum graphs in the limit of large mass,} 
Ann. Inst. H. Poincare C Anal. Non Lineaire, 38, 1295--1335. 2021.

\bibitem{BlaExn08} 
J. Blank, P. Exner and M. Havlicek. 
{\it Hilbert Space Operators in Quantum Physics,} 
2nd ed., Theor. Math. Phys., Springer, New York. 2008.





\bibitem{BurCas01} 
R. Burioni, D. Cassi, M. Rasetti, P. Sodano and A. Vezzani. 
{\it Bose-Einstein condensation on inhomogeneous complex networks,} 
J. Phys. B: At. Mol. Opt. Phys., 34, 4697--4710. 2001.

\bibitem{ByrdFriedman} 
P. F. Byrd and M. D. Friedman. {\it Handbook of Elliptic Integrals for Engineers and Scientists,} 
Springer, Berlin, Heidelberg, 2. 1971.



\bibitem{CFN} 
C. Cacciapuoti, D. Finco and D. Noja. 
{\it Topology induced bifurcations for the NLS on the tadpole graph,} 
Phys. Rev. E, 91, 013206. 2015.

\bibitem{CFN2} 
C. Cacciapuoti, D. Finco and D. Noja. 
{\it Ground state and orbital stability for the NLS equation on a general starlike graph with potentials,} 
Nonlinearity, 30(8), 3271--3303. 2017.


\bibitem{Caz} 
T. Cazenave. 
{\it Semilinear Schr\"odinger Equations,} 
Courant Lecture Notes in Mathematics, AMS. vol 10, Providence, 2003.

\bibitem{Chuiko} 
G. P. Chuiko, O. V. Dvornik, S. I. Shyian and Y. A. Baganov. 
{\it A new age-related model for blood stroke volume,} 
Comput. Biol. Med., 79, 144--148. 2016.

\bibitem{Crepeau} 
E. Cr\'epeau and M. Sorine. 
{\it A reduced model of pulsatile flow in an arterial compartment,} 
Chaos Solitons Fractals, 34(2), 594--605. 2007.


\bibitem{Ex} 
P. Exner. 
{\it Magnetoresonance on a lasso graph,} 
Found. Phys., 27, Art. 171. 1997.

\bibitem{ExS} 
P. Exner and P. $\check{S}$eba. 
{\it Free quantum motion on a branching graph,} 
Rep. Math. Phys., 28, 7--26. 1989.

\bibitem{ExSere} 
P. Exner and E. $\check{S}$ere$\check{s}$ov\'a. 
{\it Appendix resonances on a simple graph,} 
J. Phys. A, 27, 8269--8278. 1994.

\bibitem{Fid15} 
F. Fidaleo. 
{\it Harmonic analysis on inhomogeneous amenable networks and the Bose-Einstein condensation,} 
J. Stat. Phys., 160, 715--759. 2015.

\bibitem{GSS2} 
M. Grillakis, J. Shatah and W. Strauss. 
{\it Stability theory of solitary waves in the presence of symmetry, II,} 
J. Funct. Anal., 94, 308--348. 1990.

\bibitem{GSS1} 
M. Grillakis, J. Shatah and W. Strauss. 
{\it Stability theory of solitary waves in the presence of symmetry, I,} 
J. Funct. Anal., 74(1), 160--197. 1987.


\bibitem{HPW}  D. Henry, J. Perez, and W.  Wreszinski, {\it Stability theory for solitary-wave solutions of scalar field equations}, Comm. Math. Phys. 85 (1982), no. 3, pp. 351--361. 


\bibitem{KP} 
A. Kairzhan and D. E. Pelinovsky. 
{\it Multi-pulse edge-localized states on quantum graphs,} 
Anal. Math. Phys., 11, 171 (26pp). 2021.

\bibitem{KPG} 
A. Kairzhan, D. E. Pelinovsky and R. Goodman. 
{\it Drift of spectrally stable shifted states on star graphs,} 
SIAM J. Appl. Dyn. Syst., 18, 1723--1755. 2019.

\bibitem{KMPX} 
A. Kairzhan, R. Marangell, D. E. Pelinovsky and K. Xiao. 
{\it Existence of standing waves on a flower graph,} 
J. Differential Equations, 271, 719--763. 2021.

\bibitem{KNP} 
A. Kairzhan, D. Noja and D. E. Pelinovsky. 
{\it Standing waves on quantum graphs,} 
J. Phys. A: Math. Theor., 55, 243001 (51pp). 2022.

\bibitem{K} 
P. Kuchment. 
{\it Quantum graphs, I. Some basic structures,} 
Waves Random Media, 14, 107--128. 2004.

\bibitem{LP}  F. Linares and G. Ponce, \textit{Introduction to Nonlinear Dispersive Equations}, 2nd edition, Universitext, Springer, New York, (2015).



\bibitem{Mug15} 
D. Mugnolo. 
{\it Mathematical Technology of Networks,} 
Springer Proc. Math. Stat., 128, Bielefeld, December 2013. 2015.








\bibitem{Noj14}
D. Noja. 
{\it Nonlinear Schr\"odinger equation on graphs: recent results and open problems,} 
Philos. Trans. R. Soc. Lond. Ser. A Math. Phys. Eng. Sci., 372, 20130002 (20pp). 2014.

\bibitem{NPS} 
D. Noja, D. E. Pelinovsky and G. Shaikhova. 
{\it Bifurcations and stability of standing waves in the nonlinear Schr\"odinger equation on the tadpole graph,} 
Nonlinearity, 28, 2343--2378. 2015.

\bibitem{NP}  
D. Noja and D. E. Pelinovsky. 
{\it Standing waves of the quintic NLS equation on the tadpole graph,} 
Calc. Var., 59, 173. 2020.


\bibitem{Pan}  
A. Pankov. 
{\it Nonlinear Schr\"odinger equations on periodic metric graphs,} 
Discrete Contin. Dyn. Syst. A, 38(4), 697--714. 2018.

\bibitem{pazy} 
A. Pazy. 
{\it Semigroups of linear operators and applications to partial differential equations,} 
Vol. 44. Applied Mathematical Sciences, Springer-Verlag, New York, 1983.



\bibitem{Schu2015} 
C. Schubert, C. Seifert, J. Voigt and M. Waurick. 
{\it Boundary systems and (skew-) self-adjoint operators on infinite metric graphs,} 
Math. Nachr., 288, 1776--1785. 2015.

\bibitem{SBM} 
Z. A. Sobirov, D. Babajanov and D. Matrasulov. 
{\it Nonlinear standing waves on planar branched systems: Shrinking into metric graph,} 
Nanosystems, 8, 29--37. 2017.




\end{thebibliography}
\end{document}